\documentclass[10pt]{article}
\date{October 6, 2010} 

\usepackage{amssymb}
\usepackage{amsmath}
\input{liemacs.sty} 

\begin{document} 


\title{On Analytic Vectors for Unitary Representations \\ 
of Infinite Dimensional Lie Groups} 
\author{Karl-Hermann Neeb\begin{footnote}{
Department  Mathematik, FAU Erlangen-N\"urnberg, Bismarckstrasse 1 1/2, 
91054-Erlangen, Germany; neeb@mi.uni-erlangen.de}
\end{footnote}
\begin{footnote}{Supported by DFG-grant NE 413/7-1, Schwerpunktprogramm 
``Darstellungstheorie''.} 
\end{footnote}}

\maketitle
\date 
\begin{abstract} Let $G$ be a connected and simply connected 
Banach--Lie group or, more generally, 
a BCH--Lie group. On the complex enveloping algebra 
$U_\C(\g)$ of its Lie algebra $\g$ we define the concept of an 
analytic functional and show that every positive analytic 
functional $\lambda$ is integrable in the sense that it is of the form 
$\lambda(D) = \la \dd\pi(D)v, v\ra$ for an analytic vector 
$v$ of a unitary representation of $G$. 
On the way to this result we derive criteria 
for the integrability of $*$-representations of infinite 
dimensional Lie algebras of unbounded operators to unitary 
group representations. 

For the matrix coefficient $\pi^{v,v}(g) = \la \pi(g)v,v\ra$ 
of a vector $v$ in a unitary representation of an 
analytic Fr\'echet--Lie group $G$ we show that $v$ is an 
analytic vector if and only if $\pi^{v,v}$ is analytic 
in an identity neighborhood. Combining this insight with the 
results on positive analytic functionals, we derive that 
every local positive definite analytic function on a 
$1$-connected Fr\'echet--BCH--Lie group $G$ extends to a 
global analytic function. 

\end{abstract} 

\tableofcontents
\section{Introduction} 

Let $G$ be an analytic Lie group modeled on a 
locally convex space and $(\pi, \cH)$ be a unitary 
representation of $G$.\begin{footnote}{Cf.~\cite{Ne06} for a survey on locally convex 
Lie theory and Section~\ref{sec:2} for basic definitions.} 
\end{footnote}We call a vector $v \in \cH$ 
smooth, resp., analytic if its orbit map 
$\pi^v \: G \to \cH, g \mapsto \pi(g)v$ 
is smooth, resp., analytic. We write $\cH^\infty$, resp., 
$\cH^\omega$ for the subspace of smooth, resp., analytic 
vectors. Analytic vectors for unitary representations 
of finite dimensional Lie groups play an important role 
in the translation process between the representation 
$\pi$ of a connected Lie group 
$G$ and the derived representation $\dd\pi \: \g \to \End(\cH^\infty)$ 
of its Lie algebra $\g$. 
In particular, for a subspace $E \subeq \cH$ consisting of smooth 
vectors and invariant under $\g$, we can only conclude its invariance 
under $G$ 
if the set of analytic vectors in $E$ is dense 
(\cite[Cor.~4.4.5.5]{Wa72}). 
The key point in the argument is the Taylor expansion of an 
analytic function $f \: G \to E$ with values in a Banach space $E$ 
near the identity:  
$$ f(\exp X) = \sum_{n = 0}^\infty \frac{1}{n!} (X^n f)(\1), $$
where the elements $X \in \g$ acts as a left invariant 
differential operator on $C^\infty(G,E)$. The summands in this expansion 
define a linear map $T(f) \: U(\g) \to E, D \mapsto (Df)(\1)$. 
Our motivation for this paper was to understand which functionals 
on $U(\g)$ correspond in this sense to a matrix coefficient  
$\pi^{v,v}(g) = \la \pi(g)v,v\ra$ of an analytic 
vector $v$ for a unitary representation $(\pi, \cH)$ of $G$. 

As we shall see below, all this concerns only representations 
$(\pi, \cH)$ for which $\cH^\omega$ is dense. For finite dimensional 
Lie groups, this is always the case for continuous unitary 
representations on Hilbert spaces. 
Building on previous work of Harish--Chandra on semisimple groups, 
Cartier and Dixmier prove in \cite{CD58} the density 
of $\cH^\omega$ for unitary representations and 
continuous Banach representations which 
are bounded on a certain discrete central subgroup. 
The general case was obtained by 
Nelson (\cite[Thm.~4]{Nel59}) and 
G\aa{}rding (\cite{Ga60}) by convolution with heat kernels 
(cf.\ \cite[Sect.~4.4.5]{Wa72}). 
For generalizations of these density results to representations 
on locally convex spaces, we refer to \cite{Mo72}. 

For infinite dimensional Lie groups a continuous 
unitary representation need not possess any differentiable vector. 
Refining a construction from \cite{BN08}, we have shown 
in \cite{Ne10} that the continuous unitary representation 
of the abelian Banach--Lie group $G = (L^p([0,1],\R),+)$, $p \geq 1$,  
on the Hilbert  space $\cH = L^2([0,1],\C)$ by $\pi(g)f = e^{ig}f$ 
satisfies $\cH^\infty = \{0\}$. For $p = 1$, we thus obtain a continuous 
representation with no non-zero $C^1$-vector, and for 
$p = 2k$, the space of $C^k$-vectors is non-zero, but 
there is no non-zero $C^{k+1}$-vector.

As a consequence of these examples, one cannot expect 
analytic vectors to form a dense subspace for infinite dimensional 
Lie groups, and not even for Banach--Lie groups. Up to now, 
only a few classes of groups are known for which 
$\cH^\omega$ is dense for every continuous unitary representation. 
If $G$ is the additive group of a nuclear locally convex space, 
this follows from the Bochner--Minlos Theorem and 
the regularity of the corresponding measures 
(cf.\ \cite[Thm.~2.11]{Heg72} and \cite[Sect.~2.5]{Sa91} for 
the special case of the countably dimensional vector group $G = \R^{(\N)}$). 
One also has a corresponding result 
for Heisenberg groups: Let $V$ be a locally 
convex space, endowed with a continuous positive definite 
scalar product. Then we define the corresponding 
Heisenberg group by 
$$ \Heis(V) = \R \times V \times V, \quad 
(z,v,w)(z',v',w') := 
(z + z' + \shalf(\la v,w'\ra - \la v',w\ra), v + v', w + w'). $$
In \cite{Heg72} Hegerfeldt shows that, if $V$ is separable, barrelled 
and nuclear, for any continuous 
unitary representation $(\pi, \cH)$ of $\Heis(V)$, 
the space $\cH^\omega$ is dense. He even proves the existence 
of a dense subspace $\cD \subeq \cH^\omega$ 
with the property that for every $v \in \cD$ 
the orbit map $\pi^v(g) = \pi(g)v$ extends to a 
holomorphic map $\Heis(V)_\C \to \cH$ (see also \cite{Re69} for the 
special case $V = \R^{(\N)}$ and certain Banach completions). 
Further results for direct limits of finite dimensional 
matrix groups can be found in \cite[Sect.~6.4]{Sa91}. 

In this paper we are concerned with two aspects of analytic 
vectors for representations of a $1$-connected Lie group 
$G$ with Lie algebra~$\g$: 
\begin{description}
\item{(IR)} The integrability problem for Lie algebra 
representations: Given a $*$-representation 
$(\rho, \cD)$ of $\g$ on a pre-Hilbert space $\cD$, when 
is there a continuous unitary representation $(\pi, \cH)$ of $G$ 
on the completion $\cH$ of $\cD$ with 
$\cD\subeq \cH^\infty$ and 
$\rho(x) = \dd\pi(x)\res_{\cD}$ for $x \in \g$? The criteria we discuss here 
are based on a rich supply of analytic vectors. 
\item{(IF)} The integrability problem for functionals on 
$U(\g)$: Which linear functionals 
$\lambda \: U_\C(\g) := U(\g)_\C \to \C$ are {\it integrable} 
in the sense that there exists a continuous unitary 
representation $(\pi, \cH)$ of $G$ and a smooth vector 
$v \in \cH^\infty$ with $\lambda(D) = \la \dd\pi(D)v,v\ra$ 
for $D \in U(\g)$? 
\end{description}

As far as we know, the present knowledge on (IR) is limited 
to finite dimensional Lie algebras and groups. 
It is based on Nelson's Theorem, asserting that a symmetric 
operator $A$ on a Hilbert space $\cH$ is essentially selfadjoint 
if it has a dense subspace of vectors $v$ which are analytic 
in the sense that $\sum_n \frac{\|A^nv\|t^n}{n!} < \infty$ 
holds for some $t > 0$. In view of Stone's Theorem, 
for such operators $A$, 
the selfadjoint closure $\oline{A}$ generates a unitary 
one-parameter group $e^{it\oline A}$. 
The first generalization of this 
result to finite dimensional Lie algebras $\g$ can be found in 
\cite{FSSS72}, where it is shown that any representation 
$\rho \: \g \to \End(\cD)$  by skew-symmetric operators on a dense subspace 
$\cD$ of a Hilbert space $\cH$ can be integrated to a unitary group 
representation if $\cD$ consists of analytic vectors for 
the operators $\rho(x_1), \ldots, \rho(x_n)$, where 
$x_1, \ldots, x_n$ is a basis of $\g$. 
In \cite{Si72} Simon generalizes this result to the situation 
where $x_1, \ldots, x_n$ are only required to generate $\g$ as a 
Lie algebra. These integrability criteria for 
representations of finite dimensional Lie algebras 
by unbounded operators generalize to representations 
on Banach spaces (cf.\ \cite{GJ83}, see also \cite[Cor. A.4]{Jo88}). 
We recommend Section~5 in \cite{Mag92} for a detailed 
discussion of analytic vectors and integrability conditions for 
finite dimensional Lie algebras. One of our main results 
is a solution of (IR) for Banach--Lie algebras asserting the 
integrability of $\rho$ if $\cD$ consists of analytic vectors. 

Problem (IF) is an infinite dimensional variant of the 
classical (non-commutative) moment problem. 
Since smoothness of a vector $v \in \cH$ is equivalent to smoothness 
of the corresponding positive definite function 
$\pi^{v,v}(g) = \la \pi(g)v,v\ra$ (\cite[Thm.~7.2]{Ne10}), 
the GNS construction implies that 
(IF) is equivalent to the existence of a smooth positive definite 
function $\phi$ with $\lambda(D) = (D\phi)(\1)$ for $D \in U_\C(\g)$. 

If $\g$ is finite dimensional and $x_1, \ldots, x_k$ 
is a basis of $\g$, then every linear functional 
$\lambda \in U_\C(\g)^*$ defines a function 
$$ s_\lambda \: \N_0^k \to \C, \quad 
\alpha \mapsto \lambda(x^\alpha), \quad 
x^\alpha := x_1^{\alpha_1} \cdots x_k^{\alpha_k}. $$
In view of the Poincar\'e--Birkhoff--Witt Theorem, 
the correspondence between linear functionals on $U_\C(\g)$ 
and functions on $\N_0^k$ is one-to-one. The 
{\it non-commutative moment problem} is to decide 
for $s \: \N_0^k \to \C$ when the corresponding functional 
$\lambda_s \in U_\C(\g)^*$ is integrable 
in the sense of (IF) (cf. \cite[Sect.~12.2]{Sa91}). 
For $G = \R^k$ this specializes to the classical moment problem 
which, for $(s_\alpha)_{\alpha \in \N_0^k}$, asks for a 
bounded positive measure $\mu$ on $\R^k$ with 
$$ s_\alpha = \int_{\R^k} x^\alpha \dd\mu(x) \quad \mbox{ for } \quad 
\alpha \in \N_0^k.$$ 
Here we use Bochner's Theorem which says that the Fourier transform 
$\mu \mapsto \hat\mu$ defines a bijection between bounded positive 
Borel measures on $\R^k$ and continuous positive definite functions 
on $\R^k$ 
(see \cite{Ak65} for more on the classical moment problem). 
For results concerning moment problems for nuclear spaces,  
we refer to \cite[Thm.~12.5.2]{Sch90} and the original 
papers \cite{BY75}, \cite{Heg75}. 

An obvious necessary condition for the solvability of (IF) is that 
$\lambda$ must be a positive functional on the $*$-algebra 
$U_\C(\g)$ in the sense that $\lambda(D^*D) \geq 0$ for 
$D \in U_\C(\g)$, but in general this condition is not sufficient 
for integrability. For $G = \R^k$, a necessary and sufficient condition 
for integrability is non-negativity of $\lambda$ 
on positive polynomials in $U(\g) \cong \R[x_1,\ldots, x_k]$, 
which is stronger than positivity  because, for $k \geq 2$,  
not every non-negative polynomial is a sum of squares 
(\cite[Ex.~2.6.11]{Sch90}, \cite[Sects.~6.1/2]{BCR84}). 
In the light of these results, it is quite remarkable 
that, for $G = \R^k$, positivity of $\lambda$ plus estimates on the 
$s_\alpha$ which ensure convergence of the series 
$\sum_\alpha \frac{1}{\alpha!} x^\alpha s_\alpha$ in a 
zero-neighborhood of $\R^k$ imply  
unique solvability of (IF) (\cite[Lemma~12.2.1]{Sa91}, 
\cite[Thm.~II.12.7]{Sh96}).  
A non-commutative variant of this result for Heisenberg groups 
can be found in \cite[Thm.~12.2.3]{Sa91}. 

For the solvability of (IF) for a connected Lie group 
$G$, Schm\"udgen shows in \cite[Thm.~5.1]{Sch78} that 
$\lambda \in U_\C(\L(G))^*$ is integrable if and only if 
it is strongly positive in the sense that 
$\lambda(D) \geq 0$ whenever $\dd\pi(D) \geq 0$ holds for 
all unitary representations $(\pi, \cH)$ of $G$. He also shows in 
\cite[Thm.~4.1]{Sch78} that for $G \not\cong \R$ 
there exist positive functionals  which are not strongly positive. 
With respect to (IR) for finite dimensional groups,  
one has a necessary and sufficient criterion 
in terms of so-called complete positivity 
(\cite{Sch78}, \cite[Thm.~11.4.2]{Sch90}, \cite[Thm.~4.5]{Po74}). 
Here we are only concerned with 
the unique solvability of (IF) for positive 
functionals which are ``analytic'' in the sense that they 
satisfy certain estimates. 

The plan of this paper is as follows. 
After recalling the basic concepts concerning infinite dimensional 
Lie groups and analytic functions in Section~\ref{sec:2}, 
we introduce in Section~\ref{sec:3} the concept of an analytic 
linear map $\beta \: U_\C(\g) \to E$ to a Banach space $E$. 
These are linear maps for which the series 
$\sum_{n = 0}^\infty \frac{\beta(x^n)}{n!}$ converges on a 
$0$-neighborhood of $\g$. 
The main result of Section~\ref{sec:3} 
is Theorem~\ref{thm:1.5} which asserts that, for 
any continuous seminorm $p$ on $\g$ which is submultiplicative 
in the sense that $p([x,y]) \leq p(x)p(y)$ for $x,y \in \g$, 
the space of those 
analytic linear maps for which 
$\sum_n \frac{\beta(x^n)}{n!}$ converges uniformly on the 
$r$-ball with respect to $p$ is invariant under the left and 
right regular representation of $U(\g)$ on $\Hom(U(\g),E)$. 

In Section~\ref{sec:4} we turn to representations 
$\rho \: \g \to \End(V)$ of 
a Lie algebras $\g$ on a normed space $V$. We call such a representation 
{\it strongly continuous} if the orbit maps 
$\rho^v \: \g \to V, x \mapsto \rho(x)v$ are continuous 
and say that $v \in V$ is {\it analytic} if 
$\sum_n \frac{\|\rho(x)^nv\|}{n!}$ converges on some 
$0$-neighborhood of $\g$. If $\g$ is a Baire space, 
$v$ is analytic if and only if 
it is an analytic vector for all operators $\rho(x)$ separately 
(Proposition~\ref{prop:4.10}). 
In Section~\ref{sec:5} we eventually turn to unitary representations 
of an analytic Lie group $G$ with Lie algebra $\g$ by 
showing that $v \in \cH$ is analytic if and only if 
the matrix coefficient $\pi^{v,v}(g) = \la \pi(g)v,v\ra$ is analytic 
in some $\1$-neighborhood. This reduces the analyticity condition 
for $v$ to a local requirement on a scalar-valued function. 
In Section~\ref{sec:6} we proceed to the infinitesimal picture 
with the observation that for a $*$-representation 
of $\g$ on a pre-Hilbert space $\cD$, an element 
$v \in \cD$ is an analytic vector if and only if the linear functional 
$\rho^{v,v}(D) := \la \rho(D)v,v\ra$ on $U(\g)$ is analytic. 
The main results of Section~\ref{sec:6} are integrability 
results for $*$-representations on pre-Hilbert spaces. 
Theorem~\ref{thm:2.8} asserts that for 
a $1$-connected Banach--Lie group $G$ 
a $*$-representation $(\rho, \cD)$ of $\g$ is integrable 
if $\cD$ consists of analytic vectors. 
Refining this result for the larger class of BCH--Lie groups, 
we derive in particular, that 
every positive analytic functional $\lambda$ on the 
$*$-algebra $U_\C(\g)$ is integrable in the sense of (IF). 
We conclude Section~\ref{sec:6} with a variant of the integrability 
theorem on Lie algebra representations 
that does not require an a priori given Lie group $G$. 
The interest in this refinement lies in the fact that not every 
Banach--Lie algebra is the Lie algebra of some Banach--Lie group 
(\cite{EK64}). The class of BCH--Lie groups is substantially 
larger than the class of Banach--Lie groups. It contains in particular 
all (projective limits of) nilpotent Lie groups and 
groups of smooth maps on compact manifolds with values in 
Banach--Lie groups (see \cite{GN10} for more details). 
In particular, all Heisenberg groups are BCH. 

From the integrability of positive analytic functionals 
on $U_\C(\g)$, we further derive in Section~\ref{sec:7} that 
any germ of an analytic positive definite function $\phi$ 
on a $1$-connected BCH--Lie group extends uniquely to a global 
analytic positive definite function (Theorem~\ref{thm:extension}). 

An interesting issue that deserves to be pursued in the future 
is that the concept of an analytic vector makes sense for any 
representation $(\rho,V)$ of a locally convex Lie algebra $\g$ 
on a normed space $V$, even when there is no analytic Lie group $G$
with Lie algebra $\g$. This may open a door 
to applications for large classes of non-analytic Lie groups such as 
diffeomorphism groups of compact smooth manifolds. 

In his recent paper \cite{Me10}, S.~Merigon generalizes our 
Integrability Theorem~\ref{thm:5.2} for simply connected 
Banach--Lie groups in the sense that he only requires 
the existence of a direct sum decomposition 
$\g = \g_1 \oplus \cdots \oplus \g_n$ into closed subspaces, 
such that all vectors are analytic for every $x \in \g_j$, $j =1,\ldots, n$. 
We expect the results in this paper and \cite{Me10} to 
be the key tools to generalize the theory of 
analytic extension of unitary group representations 
(\cite{LM75}, \cite{Jo86, Jo87}, \cite{Ol90, Ol91}, \cite{Ne00}) to 
the context of Banach--Lie groups.  

\section{Locally convex Lie groups} \label{sec:2}

In this section we briefly recall the basic concepts related 
to infinite dimensional Lie groups and calculus on locally convex spaces. 

\begin{defn} 
(a) Let $E$ and $F$ be locally convex spaces, $U
\subeq E$ open and $f \: U \to F$ a map. Then the {\it derivative
  of $f$ at $x$ in the direction $h$} is defined as 
$$ \dd f(x)(h) := (\partial_h f)(x) := \derat0 f(x + t h) 
= \lim_{t \to 0} \frac{1}{t}(f(x+th) -f(x)) $$
whenever it exists. The function $f$ is called {\it differentiable at
  $x$} if $\dd f(x)(h)$ exists for all $h \in E$. It is called {\it
  continuously differentiable}, if it is differentiable at all
points of $U$ and 
$$ \dd f \: U \times E \to F, \quad (x,h) \mapsto \dd f(x)(h) $$
is a continuous map. Note that this implies that the maps 
$\dd f(x)$ are linear (cf.\ \cite[Lemma~2.2.14]{GN10}). 
The map $f$ is called a {\it $C^k$-map}, $k \in \N \cup \{\infty\}$, 
if it is continuous, the iterated directional derivatives 
$$ \dd^{j}f(x)(h_1,\ldots, h_j)
:= (\partial_{h_j} \cdots \partial_{h_1}f)(x) $$
exist for all integers $j \leq k$, $x \in U$ and $h_1,\ldots, h_j \in E$, 
and all maps $\dd^j f \: U \times E^j \to F$ are continuous. 
As usual, $C^\infty$-maps are called {\it smooth}. 

  (b) If $E$ and $F$ are complex locally convex spaces, then $f$ is 
called {\it complex analytic} if it is continuous and for each 
$x \in U$ there exists a $0$-neighborhood $V$ with $x + V \subeq U$ and 
continuous homogeneous polynomials $\beta_k \: E \to F$ of degree $k$ 
such that for each $h \in V$ we have 
$$ f(x+h) = \sum_{k = 0}^\infty \beta_k(h), $$
as a pointwise limit (\cite{BS71b}). 
The map $f$ is called {\it holomorphic} if it is $C^1$ 
and for each $x \in U$ the 
map $\dd f(x) \: E \to F$ is complex linear (cf.\ \cite[p.~1027]{Mil84}). 
If $F$ is sequentially complete, then $f$ is holomorphic if and only if 
it is complex analytic (\cite[Ths.~3.1, 6.4]{BS71b}). 

(c) If $E$ and $F$ are real locally convex spaces, 
then we call a map $f \: U \to F$, $U \subeq E$ open, 
{\it real analytic} or a $C^\omega$-map, 
if for each point $x \in U$ there exists an open neighborhood 
$V \subeq E_\C$ and a holomorphic map $f_\C \: V \to F_\C$ with 
$f_\C\res_{U \cap V} = f\res_{U \cap V}$  (cf.\ \cite{Mil84}). 
The advantage of this definition, which differs from the one in 
\cite{BS71b}, is that it also works nicely for non-complete spaces. 
Any analytic map is smooth, 
and the corresponding chain rule holds without any condition 
on the underlying spaces, which is the key to the definition of 
analytic manifolds (see \cite{Gl02} for details).
\end{defn}

Once the concept of a smooth function 
between open subsets of locally convex spaces is established  
(cf.\ \cite{Ne06}, \cite{Mil84}, \cite{GN10}), it is clear how to define 
a locally convex smooth manifold. 
A {\it (locally convex) Lie group} $G$ is a group equipped with a 
smooth manifold structure modeled on a locally convex space 
for which the group multiplication and the 
inversion are smooth maps. We write $\1 \in G$ for the identity element and 
$\lambda_g(x) = gx$, resp., $\rho_g(x) = xg$ for the left, resp.,
right multiplication on $G$. Then each $x \in T_\1(G)$ corresponds to
a unique left invariant vector field $x_l$ with 
$x_l(g) := T_\1(\lambda_g)x, g \in G.$
The space of left invariant vector fields is closed under the Lie
bracket of vector fields, hence inherits a Lie algebra structure. 
In this sense we obtain on $\g := T_\1(G)$ a continuous Lie bracket which
is uniquely determined by $[x,y]_l = [x_l, y_l]$ for $x,y \in \g$. 
We shall also use the functorial 
notation $\L(G) := (\g,[\cdot,\cdot])$ 
for the Lie algebra of $G$ and, accordingly, 
$\L(\phi) = T_\1(\phi)\: \L(G_1) \to \L(G_2)$ 
for the Lie algebra morphism associated to 
a morphism $\phi \: G_1 \to G_2$ of Lie groups. 
Then $\L$ defines a functor from the category 
of locally convex Lie groups to the category of 
{\it locally convex Lie algebras}, i.e., locally convex spaces, endowed with 
continuous Lie brackets. The adjoint action of $G$ on $\L(G)$ is defined by 
$\Ad(g) := \L(c_g)$, where $c_g(x) = gxg^{-1}$ is the conjugation map. 
The adjoint action is smooth and 
each $\Ad(g)$ is a topological isomorphism of $\g$. 
If $\g$ is a Fr\'echet, resp., a Banach space, then 
$G$ is called a {\it Fr\'echet-}, resp., a 
{\it Banach--Lie group}. 

A smooth map $\exp_G \: \g \to G$  is called an {\it exponential function} 
if each curve $\gamma_x(t) := \exp_G(tx)$ is a one-parameter group 
with $\gamma_x'(0)= x$. Not every infinite dimensional Lie group has an 
exponential function (\cite[Ex.~II.5.5]{Ne06}), but exponential functions 
are unique whenever they exist. 
A Lie group $G$ is said to be 
{\it locally exponential} 
if it has an exponential function for which there is an open $0$-neighborhood 
$U$ in $\L(G)$ mapped diffeomorphically by $\exp_G$ onto an 
open subset of $G$. If, in addition, $G$ is analytic and 
the exponential function is an analytic local diffeomorphism in $0$, 
then $G$ is called a {\it BCH--Lie group}. 
Then the Lie algebra $\L(G)$ is a 
{\it BCH--Lie algebra} (for Baker--Campbell--Hausdorff), i.e., there exists an open 
$0$-neighborhood $U \subeq \g$ such that for $x,y \in U$ the BCH series 
$$x * y = x + y + \frac{1}{2}[x,y] + \cdots  $$
converges and defines an analytic function 
$U \times U \to \g, (x,y) \mapsto x * y$. 
The class of BCH--Lie groups contains in particular all Banach--Lie groups 
(\cite[Prop.~IV.1.2]{Ne06}).

If $\pi \: G \to \GL(V)$ is a representation of $G$ on a locally convex 
space $V$, the exponential function 
permits us to associate to each element $x$ of the Lie algebra a 
one-parameter group $\pi_x(t) := \pi(\exp_G tx)$. 
We therefore {\bf assume} in the following that 
$G$ has an exponential function.

\section{Analytic functionals} \label{sec:3}

In the following $\g$ denotes a real locally convex Lie algebra and 
$U(\g)$ its enveloping algebra. 
The continuity of the Lie bracket on $\g$ means that, for 
every continuous seminorm $p$ on $\g$, there exists a continuous 
seminorm $\tilde p$ on $\g$ satisfying 
\begin{equation}
  \label{eq:brack-esti}
p([x,y]) \leq \tilde p(x)\tilde p(y) \quad \mbox{ for } \quad x,y \in \g. 
\end{equation}
If this relation holds with $\tilde p = p$, then 
we call $p$ is {\it submultiplicative}. 
If $\g$ is a Banach--Lie algebra, then we assume that the norm on $\g$ 
is submultiplicative 
\begin{equation}
  \label{eq:submult}
\|[x,y]\| \leq \|x\|\cdot\|y\| \quad \mbox{ for } \quad x,y \in \g, 
\end{equation}
which can always be achieved by replacing it by a suitable multiple.

In this section we introduce for each continuous submultiplicative 
seminorm $p$ on $\g$ the concept of a $p$-analytic 
linear map $\beta$ on the enveloping algebra $U(\g)$ with values in a Banach 
space~$E$ and define its radius of convergence $R_{\beta,p}$. The main 
result of this section is Theorem~\ref{thm:1.5}, asserting that 
the subspace of elements 
$\beta$ with $R_{\beta,p} \geq r$ 
is invariant under the left and right regular action of 
$U(\g)$. This applies in particular to the norm $p(x) = \|x\|$ 
of a Banach--Lie algebra. 

\begin{defn} Let $V$ be a locally convex space and $W$ be a normed space. 
We write $\Mult^n(V,W)$ for the space of continuous 
$W$-valued $n$-linear maps on $V^n$ and 
$\Sym^n(V,W)$ for the subspace of symmetric $n$-linear maps. 
We also write 
$\cP(V)$ for the set of continuous seminorms on $V$. 

(a) For an $n$-linear map $\beta \: V^n \to W$ 
and $p \in \cP(V)$, we define 
$$ \|\beta\|_p := \sup \{ \|\beta(v_1,\ldots, v_n)\| \: 
v_1, \ldots, v_n \in V, p(v_i)\leq 1\} 
\in [0,\infty]. $$
Then $\beta$ is continuous with respect to the topology on 
$V$ defined by $p$ if and only if $\|\beta\|_p < \infty$, and we have 
$$ \|\beta(v_1,\ldots, v_n)\| 
\leq \|\beta\|_p p(v_1) \cdots p(v_n) \quad \mbox{ for } \quad v_i \in V. $$

If $V$ and $W$ are Banach spaces and $p(v) = \|v\|$, 
then we thus obtain on $\Mult^n(V,W)$ the structure of a Banach space. 
In this case we simply write $\|\beta\|$ for the norm of $\beta$ 
with respect to $p$. 

(b) For $\beta \in \Mult^n(V,W)$, we write $\beta^s$ 
for the symmetrization of $\beta$: 
$$ \beta^s(x_1,\ldots, x_n) := \frac{1}{n!} \sum_{\sigma \in S_n} 
\beta(x_{\sigma(1)}, \ldots, x_{\sigma(n)}) $$ 
and observe that
$\|\beta^s\|_p \leq \|\beta\|_p$ for every $p \in \cP(V)$. 

(c) A map $f \: V \to W$ is called {\it homogeneous of degree $n$} 
if there exists a $\tilde f \in \Sym^n(V,W)$ with 
$$ f(x) = \tilde f(x,x,\ldots, x) \quad \mbox{ for }  \quad x \in V. $$
Any such $f$ is smooth and $\tilde f$ can be recovered from $f$ by 
$$ \tilde f(v_1,\ldots, v_n) 
= \frac{1}{n!}(\partial_{v_1} \cdots \partial_{v_n} f)(0). $$
For $p \in \cP(V)$, we define 
$$ \|f\|_p := \sup \{ \|f(x)\| \: p(x) \leq 1 \}. $$
Then we clearly have $\|f\|_p \leq \|\tilde f\|_p$, 
and, in view of \cite[Prop.~1.1]{BS71a}, we also have 
\begin{equation}
  \label{eq:symesti}
\|\tilde f\|_p \leq \frac{(2n)^n}{n!} \|f\|_p.
\end{equation}
\end{defn}

\begin{defn} Let $E$ be a Banach space. 

(a) We call a linear map 
$\beta \: U(\g) \to E$ {\it continuous} if, for each $n \in \N_0$,  
the $n$-linear map 
$$ \beta_n \: \g^n \to E, \quad (x_1,\ldots, x_n) \mapsto 
\beta(x_1\cdots x_n) $$
is continuous and write $\Hom(U(\g),E)_c$ for the set of all 
continuous linear maps. 
If $\g$ is Banach--Lie algebra, we have on this space a family of seminorms 
$$ p_n \: \Hom(U(\g),E)_c \to \R,\quad 
\beta \mapsto \|\beta_n\|. $$

(b) We call $\beta \in \Hom(U(\g),E)_c$ {\it analytic} if the series 
$\sum_{n = 0}^\infty \frac{\beta(x^n)}{n!}$
converges on a $0$-neighborhood of $\g$. The set of analytic 
elements is denoted $\Hom(U(\g),E)_\omega$, and for 
$E = \C$ the elements of $\Hom(U(\g),\C)_\omega$ are called 
{\it analytic functionals}. 
\end{defn} 

\begin{lem} \label{lem:frechet} {\rm(\cite[Lemma~2.3]{MNS09})} 
If $\g$ is a Banach--Lie algebra, 
then the sequence $(p_n)_{n \in \N_0}$ of 
seminorms turns the space $\Hom(U(\g),E)_c$ into a Fr\'echet space 
for which all evaluation maps 
$\ev_D \: \Hom(U(\g),E)_c \to E, \beta \mapsto \beta(D)$
are continuous. The left regular representation of $U(\g)$ on this space 
defined by 
$$ \g \times \Hom(U(\g),E)_c \to \Hom(U(\g),E)_c, \quad 
(x,\beta) \mapsto -\beta \circ \lambda_x, \quad 
\lambda_x(D) = xD, $$
and the right regular representation 
$$ \g \times \Hom(U(\g),E)_c \to \Hom(U(\g),E)_c, \quad 
(x,\beta) \mapsto \beta \circ \rho_x, \quad 
\rho_x(D) = Dx, $$
are continuous bilinear maps. 
\end{lem} 

\begin{prop} \label{prop:anachar} 
If $\g$ is a Baire space, f.i., if it is Fr\'echet, 
then an element $\beta \in \Hom(U(\g),E)_c$ is analytic 
if and only if there exists a $p \in \cP(\g)$ 
with 
\begin{equation}
  \label{eq:convsa}
 \sum_{n = 0}^\infty \frac{\|\beta_n^s\|_p}{n!}< \infty. 
\end{equation}
If $\g$ is Banach, this is equivalent to the existence of a 
$t > 0$ with 
\begin{equation}
  \label{eq:convs}
 \sum_{n = 0}^\infty \frac{\|\beta_n^s\|t^n}{n!}< \infty. 
\end{equation}
\end{prop}

\begin{prf} From \eqref{eq:convsa} we obtain for $p(x) < 1$ the estimate 
$$  \sum_{n = 0}^\infty \frac{\|\beta(x^n)\|}{n!}
=  \sum_{n = 0}^\infty \frac{\|\beta_n^s(x,\ldots, x)\|}{n!}
\leq  \sum_{n = 0}^\infty \frac{\|\beta_n^s\|_p p(x)^n}{n!} 
\leq  \sum_{n = 0}^\infty \frac{\|\beta_n^s\|_p}{n!} < \infty. $$
Therefore \eqref{eq:convsa} implies the uniform convergence 
of the series $\sum_n \frac{\beta(x^n)}{n!}$ for $p(x)< 1$, 
hence in particular the analyticity of $\beta$. 
If $\g$ is Banach, we obtain from \eqref{eq:convs} 
the uniform convergence for $\|x\| < t$. 

Suppose, conversely, that $\beta$ is analytic. 
Then the functions $f_n(x) := \frac{\beta(x^n)}{n!}$ on $\g$ 
are continuous and homogeneous of degree~$n$. 
In view of \cite[Prop.~5.2]{BS71b}, the convergence of the sequence 
$\sum_n f_n(x)$ on some $0$-neighborhood of $\g$ implies the existence 
of some $q \in \cP(\g)$ for which the series converges 
uniformly for $q(x) \leq 1$. In particular, there exists a $C > 0$ with 
$\|\beta(x^n)\| \leq C n!$ for $q(x) \leq 1$. 
We thus obtain with 
\eqref{eq:symesti} the relation 
$$ \|\beta_n^s\|_q 
\leq \frac{(2n)^n}{n!} \|n! f_n\|_q
\leq n^n 2^n \|f_n\|_q \leq n^n 2^n C. $$
Stierling's Formula 
\begin{equation}
  \label{eq:stierling}
   n! \approx \sqrt{2 \pi n} \; \left(\frac{n}{e}\right)^{n}, 
\end{equation}
leads to 
$\frac{n^n}{n!} \approx \frac{e^n}{\sqrt{2 \pi n}},$ 
so that we find with some $C' > 0$ 
$$ \sum_{n = 0}^\infty \frac{\|\beta_n^s\|_q t^n}{n!} 
\leq C' \sum_{n = 0}^\infty \frac{e^n}{\sqrt{n}} (2t)^n 
= C' \sum_{n = 0}^\infty \frac{(2et)^n}{\sqrt{n}}. $$
This series converges for $t < \frac{1}{2e}$, so that the assertion 
follows for any such $t$ with $p := t^{-1}q$. 
\end{prf} 

\begin{defn} If $\beta \in \Hom(U(\g),E)_c$ is analytic 
and $\sum_{n = 0}^\infty \frac{\|\beta_n^s\|_p t^n}{n!}< \infty$ for some 
$t > 0$, then 
$$\sum_{n = 0}^\infty \frac{\|\beta(x^n)\|}{n!}
\leq \sum_{n = 0}^\infty \frac{\|\beta_n^s\|_p p(x)^n}{n!}< \infty $$
holds for $p(x) < t$ and 
$$ R_{\beta,p} := \sup \Big\{ t > 0 \: 
\sum_{n = 0}^\infty \frac{\|\beta_n^s\|_p t^n}{n!}< \infty\Big\} 
\in ]0,\infty] $$
is called the {\it $p$-radius of convergence of $\beta$}. 
If $R_{\beta, p}< \infty$, we say that $\beta$ is {\it $p$-analytic}, or, 
{\it analytic with respect to~$p$}. According 
to Hadamard's Formula, we have
$$ R_{\beta,p}^{-1} = \limsup \root n \of {\frac{\|\beta_n^s\|_p}{n!}}. $$
Note that, for each $r > 0$, the set 
$\Hom(U(\g),E)_{\omega,p,r}$ of analytic maps $\beta$ 
with $R_{\beta,p} \geq r$ is a linear subspace. 
We write 
$\Hom(U(\g),E)_{\omega,p}$ for the union of all these spaces, 
which also is a linear subspace. 
\end{defn}

We now come to the main result of this section. 

\begin{thm} \label{thm:1.5} 
If $p \in \cP(\g)$ is submultiplicative and $r > 0$, 
then $\Hom(U(\g),E)_{\omega,p,r}$ 
is invariant under the left and right regular action of $\g$. 
In particular, the space $\Hom(U(\g),E)_{\omega,p}$ of analytic elements 
is invariant. 
\end{thm} 

\begin{prf} For $\beta \in \Hom(U(\g), E)_{\omega,p,r}$, 
$y \in \g$, $n \in \N_0$, and 
$1 \leq k \leq n+1$, we define 
$$ (i_y^k \beta)_n\in \Mult^n(\g,E), \quad 
(i_y^k\beta)_n(x_1,\ldots, x_n) 
:= \beta(x_1\cdots  x_{k-1} y x_k \cdots x_n). $$
For $n \in \N_0$ and $y \in \g$ we put 
$$ c_n := \inf \{ C > 0 \: (\forall y \in \g)(\forall k \leq n+1)\ 
\|(i_y^k \beta)_n^s\|_p \leq C p(y)\}. $$

Now let $y \in \g$. Our first goal is an estimate for the numbers $c_n$. 
For $n = 0$ we have $k = 1$ and 
$i_y^1 \beta= \beta(y)$, so that 
$$ c_0 = \|\beta_1\|_p = \|\beta_1^s\|_p. $$
Our strategy is to compare the two symmetric $n$-linear maps 
$(i_y^k\beta)_n^s$ and $i_y(\beta_{n+1}^s) = \beta_{n+1}^s(y,\cdots)$:
\begin{align*}
& (i_y^k\beta)_{n}^s(x_1, \ldots, x_n) - 
\beta_{n+1}^s(y,x_1, \ldots, x_n)   \\ 
&= \frac{1}{n!}\Big(\sum_{\sigma \in S_n} \beta(x_{\sigma(1)} \cdots 
x_{\sigma(k-1)} y x_{\sigma(k)} \cdots x_{\sigma(n)}) \\
&\qquad - \frac{1}{n+1}\sum_{j = 1}^{n+1} \sum_{\sigma \in S_n} 
\beta(x_{\sigma(1)} \cdots x_{\sigma(j-1)} y x_{\sigma(j)} 
\cdots x_{\sigma(n)})\Big) \\ 
&= \frac{1}{(n+1)!}\sum_{j = 1}^{n+1} \sum_{\sigma \in S_n} \\ 
&\qquad \Big(\beta(x_{\sigma(1)} \cdots 
x_{\sigma(k-1)} y x_{\sigma(k)} \cdots x_{\sigma(n)}) 
- \beta(x_{\sigma(1)} \cdots x_{\sigma(j-1)} y x_{\sigma(j)} 
\cdots x_{\sigma(n)})\Big).
\end{align*}
With the relation 
$[y,x_1 \cdots x_a] = \sum_{j = 1}^a x_1 \cdots x_{j-1} [y,x_j] 
x_{j+1} \cdots x_a$
we obtain for $k < j$ 
\begin{align*}
& \beta(x_{\sigma(1)} \cdots 
x_{\sigma(k-1)} y x_{\sigma(k)} \cdots x_{\sigma(n)}) 
- \beta(x_{\sigma(1)} \cdots x_{\sigma(j-1)} y x_{\sigma(j)} 
\cdots x_{\sigma(n)})\\
&= \beta(x_{\sigma(1)} \cdots 
x_{\sigma(k-1)} [y, x_{\sigma(k)} \cdots x_{\sigma(j-1)}] 
x_{\sigma(j)} \cdots x_{\sigma(n)}) \\
&= \sum_{\ell = k}^{j-1} \beta(x_{\sigma(1)} \cdots 
x_{\sigma(\ell-1)} [y, x_{\sigma(\ell)}]x_{\sigma(\ell+1)} 
\cdots x_{\sigma(n)}).  
\end{align*}
Next we observe that 
\begin{align*}
& \sum_{\sigma \in S_n} \beta(x_{\sigma(1)} \cdots 
x_{\sigma(\ell-1)} [y, x_{\sigma(\ell)}]x_{\sigma(\ell+1)} 
\cdots x_{\sigma(n)})\\
&=  (n-1)!\sum_{m = 1}^n (i_{[y,x_m]}^\ell \beta)_{n-1}^s(x_1, \ldots, 
\hat{x_m}, \ldots, x_n),
\end{align*}
where $\hat{x_m}$ denotes omission of $x_m$. 
The norm of this sum can be estimated by 
\begin{align*}
& (n-1)! \sum_{m = 1}^n c_{n-1} p([y,x_m]) p(x_1) \cdots p(\hat x_m) 
\cdots p(x_n) \\
&\leq (n-1)! \sum_{m = 1}^n c_{n-1} p(y)p(x_1) \cdots p(x_n) 
= n! c_{n-1} p(y) p(x_1) \cdots p(x_n). 
\end{align*}
For $k > j$ we likewise obtain 
\begin{align*}
& \beta(x_{\sigma(1)} \cdots 
x_{\sigma(k-1)} y x_{\sigma(k)} \cdots x_{\sigma(n)}) 
- \beta(x_{\sigma(1)} \cdots x_{\sigma(j-1)} y x_{\sigma(j)} 
\cdots x_{\sigma(n)})\\
&= \beta(x_{\sigma(1)} \cdots 
x_{\sigma(j-1)} [x_{\sigma(j)} \cdots x_{\sigma(k-1)},y] 
x_{\sigma(k)} \cdots x_{\sigma(n)}) \\
&= \sum_{\ell = j}^{k-1} \beta(x_{\sigma(1)} \cdots 
x_{\sigma(\ell-1)} [x_{\sigma(\ell)},y]x_{\sigma(\ell+1)} 
\cdots x_{\sigma(n)})  
\end{align*}
which can be estimated in the same way. 
This leads to the estimate 
\begin{align*}
& \|(i_y^k\beta)_{n}^s(x_1, \ldots, x_n) - 
\beta_{n+1}^s(y,x_1, \ldots, x_n)\|   \\ 
&\leq \frac{1}{(n+1)!}\Big(\sum_{j > k} \sum_{\ell = k}^{j-1} 
n! c_{n-1} p(y)p(x_1) \cdots  p(x_n) + \sum_{j < k} \sum_{\ell = j}^{k-1} 
n! c_{n-1} p(y) p(x_1) \cdots p(x_n) \Big) \\ 
&\leq \frac{1}{n+1}\sum_{j \not= k} |k-j| c_{n-1}
p(y) p(x_1) \cdots p(x_n) 
\leq \frac{n^2}{n+1}c_{n-1} p(y) p(x_1) \cdots p(x_n)\\
&\leq n c_{n-1} p(y) p(x_1) \cdots p(x_n). 
\end{align*}
We conclude that 
$$ \|(i_y^k\beta)_{n}^s(x_1, \ldots, x_n)\|
\leq \|\beta_{n+1}^s\|_p p(y) p(x_1) \cdots p(x_n)  
+ n c_{n-1}  p(y) p(x_1) \cdots  p(x_n), $$
i.e., 
\begin{equation}
  \label{eq:recu1}
c_{n} \leq  \|\beta_{n+1}^s\|_p +  n c_{n-1}. 
\end{equation}
Iterating this estimate leads to 
\begin{align*}
c_n &\leq  \|\beta_{n+1}^s\|_p 
+  n \|\beta_{n}^s\|_p 
+  n(n-1)\|\beta_{n-1}^s\|_p\\
&\qquad \qquad 
+ \cdots + n(n-1)\cdots 2 \cdot \|\beta_2^s\|_p + n!\|\beta_1^s\|_p. 
\end{align*}

For $0 < t < \min(1,R_{\beta,p})$, there exists a $C> 0$ with 
$$ \|\beta_n^s\|_p \leq C n! t^{-n} \quad \mbox{ for } \quad n \in \N. $$
We thus obtain
\begin{align*}
c_n  
&\leq  C\big(
(n+1)! t^{-n-1} +  n! n t^{-n} 
+ \cdots + n\cdots 2 \cdot 2! t^{-2}  + n! t^{-1}\big) \\
&\leq  C (n+1) \cdot (n+1)! t^{-n-1} \leq C (n+2)! t^{-n-1}. 
\end{align*}
This implies that 
$$ \sum_n \frac{c_n s^n}{n!}
\leq C \sum_n (n+2)(n+1)t^{-1} (s/t)^n, $$
which converges for $s < t$. 

Finally we observe that 
$(\beta \circ \rho_y)_n^s = (i_y^{n+1}\beta)_n^s,$ 
so that 
$$ \|(\beta \circ \rho_y)_n^s\|_p = \|(i_y^{n+1}\beta)_n^s\|_p  
\leq c_n p(y) $$
shows that $\beta \circ \rho_y$ is analytic 
with $R_{\beta \circ \rho_y,p} \geq R_{\beta,p}$. 
Similarly $(\beta \circ \lambda_y)_n^s = (i_y^1 \beta)_n^s$ 
implies that $\beta \circ \lambda_y$ is analytic 
with $R_{\beta \circ \lambda_y,p} \geq R_{\beta,p}$. 
\end{prf}

\begin{rem} It is an interesting question whether (for Banach--Lie 
algebras) the condition 
$\sum_{n = 0}^\infty \frac{\|\beta_n^s\| t^n}{n!}< \infty$
characterizing analytic elements in $\Hom(U(\g),E)_c$ 
also implies the stronger condition 
$\sum_{n = 0}^\infty \frac{\|\beta_n\| t^n}{n!}< \infty$ 
on the non-symmetrized maps $\beta_n$. 
A natural way to verify this would be to compare the norms of 
$\beta_n$ and $\beta_n^s$. 
For $x_1, \ldots, x_n \in \g$ we have 
$$ (\beta_n - \beta_n^s)(x_1,\ldots, x_n) 
= \frac{1}{n!} \sum_{\sigma \in S_n} 
\beta(x_1 \cdots x_n) - \beta(x_{\sigma(1)} \cdots x_{\sigma(n)}). $$
Writing a permutation in terms of transpositions 
$s_j$ exchanging $j$ and $j+1$ and using that the corresponding 
word length satisfies $\ell(\sigma) \leq \frac{n(n-1)}{2}$, we 
obtain as in the proof of Theorem~\ref{thm:1.5} an estimate of the form 
$$ \|\beta_n - \beta_n^s\| \leq \frac{n(n-1)}{2} \|\beta_{n-1}\|, $$
which leads to the recursive estimate 
$$ \|\beta_n\| \leq \|\beta_n^s\| + \frac{n(n-1)}{2} \|\beta_{n-1}\|. $$
Iteration now yields 
\begin{align*}
&\|\beta_n\| \\
&\leq \|\beta_n^s\| 
+ \frac{n(n-1)}{2} \|\beta_{n-1}^s\|
+ \frac{n(n-1)^2(n-2)}{2^2} \|\beta_{n-2}^s\|
+ \cdots  
+ \frac{n(n-1)^2\cdots 2^2\cdot 1)}{2^{n-1}} \|\beta_1\| \\
&\leq \sum_{k = 0}^{n-1} \frac{(n!)^2}{((n-k)!)^2 2^k}\|\beta_{n-k}^s\|. 
\end{align*}
From the estimate $\|\beta_n^s\| \leq C n! t^{-n}$, we thus obtain 
\begin{align*}
\sum_n \frac{\|\beta_n\| t^n}{n!} 
&\leq C\sum_n \frac{t^n}{n!} 
\sum_{k = 0}^{n-1} \frac{(n!)^2}{((n-k)!)^2 2^k} t^{k-n} (n-k)! \\
&\leq C\sum_n \sum_{k = 0}^{n-1} \frac{n!}{(n-k)!} (t/2)^k 
= C\sum_k (t/2)^k \sum_{n > k} \frac{n!}{(n-k)!} = \infty.
\end{align*}
Apparently this estimate does not lead to the desired convergence, 
and we do not know if the series 
$\sum_n \frac{\|\beta_n\| t^n}{n!}$ converges for some $t > 0$. 

For finite dimensional Lie algebras, 
analytic vectors are usually defined by a condition equivalent to 
the convergence of the series 
$\sum_n \frac{\|\beta_n\|t^n}{n!}$ for some $t > 0$ 
(cf.\ \cite[Sect.~10.4]{Sch90}). At first sight this 
is stronger than the present one 
(cf.~\cite[Lemma~10.49]{Sch90}), but 
according to \cite[Lemma~7.1, Thm.~2]{Nel59} it is equivalent. 
We do not know if this is also true for Banach--Lie algebras. 
\end{rem}

\begin{rem} (a) If the topology on $\g$ is defined by a set of 
submultiplicative seminorms, then Theorem~\ref{thm:1.5} implies 
the invariance of the space $\Hom(U(\g),E)_\omega$ of analytic 
elements under the left and right regular action. 

(b) We do not know if this is true for general locally convex 
Lie algebras. Here the critical point is that, if 
$c_n^p$ denotes the constant $c_n$, defined with respect to the 
seminorm $p$, and 
$p([x,y]) \leq q(x)q(y)$ for some $q \in \cP(\g)$ with $q \geq p$, 
then the argument in the proof 
of Theorem~\ref{thm:1.5} only yields an estimate of the form 
$$ c_n^q \leq \|\beta_{n+1}^s\|_p+ n c_{n-1}^p. $$
We do not see how to use this estimate to conclude the 
invariance of $\Hom(U(\g),E)_{\omega}$. 
\end{rem}

\section{Analytic vectors and continuity of Lie algebra representations} \label{sec:4}

In this section $\K \in \{\R,\C\}$ and all vector spaces are 
$\K$-vector spaces. Before we turn to the connection between 
analytic functionals and unitary representations of Lie groups, 
we have to clarify some continuity requirements for 
representations of topological Lie algebras on topological 
vector spaces. 

\subsection*{Continuity of Lie algebra representations}

\begin{defn}
  \label{def:5.1}
Let $(\rho, V)$ be a representation of the topological Lie algebra 
$\g$ on the topological vector space $V$. 

(a) We say that 
$\rho$ is {\it strongly continuous} if for every $v \in V$ the map 
$$\beta_v \: \g \to V, \quad x \mapsto \rho(x)v$$ 
is continuous. For a finite dimensional Lie algebra $\g$, 
this condition is always satisfied. 

(b) If $E \subeq V'$ is a subspace separating the points of the 
completion $\hat V$ of $V$, then 
$\rho$ is called {\it weakly continuous with respect to $E$} 
if, for $\alpha \in E$ and $v \in V$, the linear functional 
$x \mapsto \alpha(\rho(x)v)$ is 
continuous.
\begin{footnote}
  {If $V$ is a normed space which is not complete and 
$\hat v \in \hat V \setminus V$, then its annihilator 
$E := \hat v^\bot \subeq V' \cong \hat V'$ is a hyperplane 
with $E^\bot \cap V = \K \hat v  \cap V = \{0\}$. 
Therefore $E$ separates the points of $V$ but not of~ $\hat V$.}
\end{footnote}

(c) If $V$ is a complex pre-Hilbert space, 
then we call $\rho$ 
\begin{description}
\item[\rm(i)] {\it weakly continuous} if, for any two vectors 
$v,w \in V$, the map $\g \to \C, x \mapsto \la \rho(x)v,w\ra$ 
is continuous, i.e., if $\rho$ is weakly continuous 
with respect to $E = \{ \la \cdot, v \ra \: v \in V\}$. 
\item[\rm(ii)] a {\it $*$-representation} if 
$$\la \rho(x)v,w \ra = - \la v, \rho(x)w\ra \quad \mbox{ for } \quad 
x \in \g, v, w \in V. $$
\end{description}
\end{defn}

\begin{lem}
 \label{lem:5.0} 
If $\g$ is a Fr\'echet space and $V$ is a metrizable locally convex 
space, then every weakly continuous representation is strongly continuous. 
\end{lem} 

\begin{prf} We claim that the graph of $\beta_v$, considered 
as a linear map into the Fr\'echet completion $\hat V$ of $V$, 
is closed. Since every continuous linear functional on $V$ extends 
uniquely to a continuous linear functional on the completion 
$\hat V$, we identify the two 
spaces $\hat V'$ and $V'$. 

Suppose that $x_n \to x$ in $\g$ and $\beta_v(x_n) \to w \in \hat V$ 
and let $E \subeq \hat V'$ be a subspace separating the points of $\hat V$ 
such that all linear functionals $\g \to \K, x \mapsto \alpha(\rho(x)v)$, 
$\alpha \in E$, $v \in V$, are continuous. 
Then $\alpha(\rho(x_n)v) \to \alpha(\rho(x)v)$ and 
$\alpha(\rho(x_n)v) \to \alpha(w)$ implies 
$w = \rho(x)v$. Therefore the graph of $\beta_v$ is closed, 
so that the Closed Graph Theorem (\cite[Ch.~1, \S 3, Cor. 5]{Bou07}) 
implies that $\beta_v$ is continuous. 
\end{prf}

The following proposition is classical for $k = 2$, but here we need a 
more general version. 

\begin{prop} \label{prop:simcont} 
If $F_1, \ldots, F_k$ are Fr\'echet spaces, $V$ a topological 
vector space, and 
$\beta \: F_1 \times\cdots \times F_k \to V$ a $k$-linear map which 
is separately continuous, then $\beta$ is continuous. 
\end{prop}

\begin{prf} We argue by induction over $k$. The 
case $k = 1$ is trivial, so that we may assume that 
$k > 1$. Let $F := F_1 \times \cdots \times F_{k-1}$ and consider the 
map 
$$ \Omega \: F \times F_k \to V, \quad 
(x,y) = ((x_1,\ldots, x_{k-1}),y)
\mapsto \beta(x_1, \ldots, x_{k-1}, y).$$
For $y \in F_k$, our induction hypothesis 
implies that the $(k-1)$-linear map 
$\Omega(\cdot, y) \: F \to V$ is continuous, so that 
the continuity of $\Omega$ follows from \cite[Prop.~5.1]{Ne10}. 
\end{prf}

\begin{lem}
  \label{lem:5.1} 
If the representation $(\rho,V)$ is strongly continuous, then, 
for each $n \in \N$ and every $v \in V$, the $n$-linear map 
$$ \beta_v^n \: \g^n \to V, \quad 
(x_1,\ldots, x_n) \mapsto \rho(x_1) \cdots \rho(x_n)v $$ 
is separately continuous. 
If, in addition, $\g$ is Fr\'echet, then each $\beta_v^n$ is continuous. 
\end{lem}

\begin{prf} We argue by induction. By definition, the assertion holds for 
$n = 1$. We assume that each map $\beta_w^{n-1}$ is 
continuous in each argument. For $j < n$, the continuity of 
$\beta_v^n$ in $x_j$ then follows from the continuity 
of $\beta_{\rho(x_n)v}^{n-1}$ in each  argument. Next we observe that 
\begin{align*}
\beta_v^n(x_1,\ldots, x_n)
&=  \beta_v^n(x_1,\ldots, x_{n-2}, x_n ,x_{n-1}) 
+ \rho(x_1) \cdots \rho(x_{n-2})\rho([x_{n-1}, x_n]) v \\
&= \beta_v^n(x_1,\ldots, x_{n-2}, x_n ,x_{n-1}) + 
\beta_v^{n-1}(x_1, \ldots, x_{n-2}, [x_{n-1}, x_n]). 
\end{align*}
From the continuity of the Lie bracket on $\g$ it now follows 
that both terms on the right are continuous in $x_n$. This 
proves the continuity of $\beta_v^n$ in $x_n$. 

Now $\beta_v^n \: \g \to V$ is a separately continuous 
multilinear map, hence continuous if $\g$ is a Fr\'echet space 
(Proposition~\ref{prop:simcont}). 
\end{prf}

\begin{lem} \label{lem:4.5} 
Suppose that $V$ is a pre-Hilbert space. 
Then a $*$-representation $(\rho,V)$ of the Fr\'echet--Lie 
algebra $\g$ on $V$ is strongly continuous if and only if 
for each $v \in V$ the quadratic map 
$q_v \: \g \to \R, x \mapsto \la \rho(x)^2v,v\ra$ is continuous.
\end{lem} 

\begin{prf} If each map $q_v$ is continuous, then 
$\|\rho(x)v\|^2 = - q_v(x)$ implies that $\rho$ is strongly continuous. 

If, conversely, $\rho$ is strongly continuous, then 
Lemma~\ref{lem:5.1} implies that the bilinear maps 
$\beta_v^2 \: \g^2 \to V, (x,y) \mapsto \rho(x)\rho(y)v$ 
are continuous, and this implies the continuity of the maps~$q_v$. 
\end{prf}

\begin{prop} If $(\pi, V)$ is a continuous representation 
of the Lie group $G$ on $V$, then the derived representation 
$\dd\pi$ of 
$\g$ on the space $V^\infty$ of smooth vectors is strongly continuous 
with respect to the subspace topology induced from~$V$. 
\end{prop}

\begin{prf} If $v$ is a smooth vector  and 
$\pi^v \: G \to V, g \mapsto \pi(g)v$, its smooth orbit map, then 
$x\mapsto \dd\pi(x)v = T_{(\1,0)}(\pi^v)(x,0)$
is a continuous linear map. 
\end{prf}

In \cite[Thm.~8.2]{Ne10} we show that the preceding observation can even 
be sharpened if $G$ and $V$ are Banach. If $\cD_\g \subeq V$ 
denotes the common domain of the generators 
$\oline{\dd\pi(x)}$ of the one-parameter groups 
$t \mapsto \pi(\exp tx)$, $x \in \g$, then, for each 
$v \in \cD_\g$, the map $\beta_v \: \g \to V$ is continuous 
and linear. In this sense the continuity requirement is satisfied 
for any representation of a Banach--Lie algebra 
that comes from a continuous Lie group action on a Banach space. 

\subsection*{Analytic vectors} 

We are now ready to turn to the central concept of this article. 

\begin{defn} Let $A$ be an operator with domain $\cD(A)$ on the 
Banach space $E$. An element in the space 
$\cD^\infty(A) := \bigcap_{n = 1}^\infty \cD(A^n)$ is called a 
{\it smooth vector for $A$}, and a smooth vector is said to be 
{\it analytic} if there exists a $t > 0$ with 
$$ \sum_{n = 0}^\infty \frac{\|A^n v\|t^n}{n!}< \infty. $$
We write $\cD^\omega(A)$ for the space of analytic vectors for~$A$. 
\end{defn}

An important application of the concept of an analytic vector is that it 
provides the following criterion for essential selfadjointness: 
\begin{thm} \label{thm:nelson} {\rm(Nelson's Theorem; \cite{Nel59}; 
\cite[Thm.~X.39]{RS75})} 
A symmetric operator $A$ on a Hilbert space $\cH$ for which 
the space $\cD^\omega(A)$ of analytic vectors is dense is essentially 
selfadjoint. 
\end{thm}

\begin{defn} Let $(\rho, V)$ be a strongly continuous 
representation of $\g$ on the normed space $V$. 
Then an element $v \in V$ is called an {\it analytic vector} if 
$\sum_{n = 0}^\infty \frac{\|\rho(x)^n v\|}{n!}$ 
converges for every $x$ in a $0$-neighborhood of $\g$. 
\end{defn}

\begin{prop} {\rm(Automatic simultaneous analyticity)} \label{prop:4.10} 
Let $(\rho, V)$ be a strongly continuous representation of the 
Lie algebra $\g$ on the normed space $V$. If $\g$ is a Baire space 
(which is in particular the case if $\g$ is Fr\'echet), then every 
$v \in V$ which is an analytic vector for all the operators 
$\rho(x)$, $x \in \g$, is an analytic vector.   
\end{prop}

\begin{prf} In view of Lemma~\ref{lem:5.1}, the functions 
$f_n(x) := \frac{\rho(x)^nv}{n!}$ are continuous and homogeneous 
of degree $n$. The assumption that $v$ is analytic for every 
operator $\rho(x)$ implies that the set $S$ of all elements 
$x \in \g$ for which the series $\sum_n f_n(x)$ converges in 
the completion $\hat V$ of $V$ is absorbing, i.e., 
$\bigcup_{n \in \N} n \cdot S = \g$. Now \cite[Prop.~5.2(1)]{BS71b} 
implies the uniform convergence of this series on a 
$0$-neighborhood of $\g$, i.e., $v$ is an analytic vector.
\end{prf}

\section{Analytic vectors of unitary representations} \label{sec:5}

The main result of this section asserts that a  vector 
$v$ in a Hilbert space $\cH$ is analytic with respect to a 
continuous unitary representation $(\pi, \cH)$ of an analytic 
Lie group~$G$ if and only if the matrix coefficient 
$\pi^{v,v}(g) = \la \pi(g)v,v\ra$ is analytic in a $\1$-neighborhood. 
This will be a direct consequence of the following theorem 
characterizing analyticity of $\cH$-valued maps in terms of 
analyticity of the corresponding kernel. 
The case of this theorem, where $M$ is one-dimensional, 
can already be found in \cite[p.~81]{Kr49}.

\begin{thm} \label{thm:kerana} Let $M$ be an analytic Fr\'echet manifold and 
$\cH$ be a Hilbert space. Then a function $\gamma \: M \to \cH$ is 
analytic if and only if the kernel 
$K(x,y) := \la \gamma(x), \gamma(y) \ra$ is analytic on an open 
neighborhood of the diagonal in $M \times M$. 
\end{thm}

\begin{prf} {\bf Step 1:} We assume w.l.o.g.\ that $\cH$ is a real 
Hilbert space. Since compositions of analytic maps are analytic 
(\cite[Prop.~2.8]{Gl02}), 
the analyticity of $\gamma$ implies the analyticity of the kernel~$K$. 
Suppose, conversely, that $K$ is analytic. As our assertion 
is local, we may assume that $M$ is an open subset 
of a locally convex space $V$ containing $0$, that 
$K$ is analytic on $M \times M$, and that we have 
to show the analyticity of $\gamma$ in a $0$-neighborhood. 
To this end, we may further assume that $U$ is a convex symmetric 
$0$-neighborhood for which $K$ extends to a holomorphic function 
on $W := (U + i U) \times (U + i U) \subeq V_\C \times V_\C$. 
Shrinking $W$ if necessary, we may further assume that 
$$C := \sup \{ |K(z,w)| \: (z,w) \in W \} < \infty. $$
\cite[Thm~7.1]{Ne10} implies that $\gamma$ is smooth. 

{\bf Step 2:} Next we consider the special case $\cH = \R$.  
Since $K$ is analytic, it is in particular smooth, so that 
$\gamma$ is a smooth curve in $\cH$, and 
$$ \|\gamma^{(n)}(t)\|^2 
= \frac{d^n}{dt^n}\frac{d^n}{ds^n}\la \gamma(t), \gamma(s)\ra\vert_{t=s} 
= (\partial_1^n \partial_2^n K)(t,t). $$
For any $t \in U$ and $r > 0$ for which 
$$ W_r(t) := \{ z \in \C^2 \: |z_1-t|, |z_2-t| \leq r \} \subeq W $$
we have the Cauchy estimates 
$$  |(\partial_1^n \partial_2^m K)(t,t)| \leq C n!m! r^{-n-m}. $$
We conclude that 
\begin{equation}
  \label{eq:derivesti}
\|\gamma^{(n)}(t)\| \leq \sqrt{C} n! r^{-n}. 
\end{equation}

{\bf Step 3:} For a general $\cH$, we consider for $x \in U$ and $h 
\in V$ the curves 
$$ \gamma_{x,h}(t) := \gamma(x + th) $$ 
which are analytic in a neighborhood of $t = 0$ by Step $2$. 
To prove the convergence of the Taylor series of $\gamma$ in a 
$0$-neighborhood, we observe that 
$$ \gamma_{x,h}^{(n)}(0) = (\dd^n\gamma)(x)(h,\ldots, h), $$
so that \eqref{eq:derivesti} implies that 
$$ \|(\dd^n\gamma)(x)(h,\ldots, h)\| 
=  \|\gamma_{x,h}^{(n)}(0)\| \leq n!2^{-n}\sqrt{C} $$
whenever $(x + z_1 h, x + z_2 h) \in W$ holds for 
$|z_i| \leq 2$. We conclude that the Taylor series 
$$\sum_{n = 0}^\infty \frac{1}{n!} (\dd^n \gamma)(x)(h,\ldots, h) $$ 
of $\gamma$ in $x$ converges uniformly for all pairs $(x,h)$ in an open 
neighborhood of $(0,0)$. 
Since the remainder term in the Taylor expansion satisfies 
\begin{align*}
& \|\gamma(x + h) - \sum_{k = 0}^n \frac{1}{n!} (\dd^n\gamma)(x)(h,\ldots, h)\|\\
&\leq \frac{1}{(n+1)!} \int_0^1 \|(\dd^{n+1} f)(x + th)(h,\ldots, h)\|\, dt
\leq 2^{-n-1} \sqrt{C} \to 0, 
\end{align*}
the Taylor series actually converges to $\gamma(x+h)$ 
for $x,h$ close to $0$. 
In view of \cite[Thm.~7.2]{BS71b}, this implies that 
$\gamma$ is analytic in the sense defined above. 
\end{prf}

The following theorem is an analytic version of 
\cite[Thm.~7.2]{Ne10}, which is concerned with smooth vectors 
(cf.\ \cite[Prop.~4.1]{Go69} for the case $G = \R$). 

\begin{thm} \label{thm:5.2} If 
$(\pi, \cH)$ is a unitary representation of an analytic Fr\'echet--Lie group 
$G$, then $v \in \cH$ is an analytic vector if and only if the 
corresponding matrix coefficient $\pi^{v,v}(g) := \la \pi(g)v,v\ra$ 
is analytic on a $\1$-neighborhood in $G$. 
\end{thm}

\begin{prf} (cf.\ \cite[Thm.~7.2]{Ne10} for smooth vectors) 
Clearly, $\pi^{v,v}$ is analytic for $v \in \cH^\omega$. 
Suppose, conversely, that $\pi^{v,v}$ is analytic in a $\1$-neighborhood 
$U \subeq G$ and that $U'$ is a $\1$-neighborhood with 
$h^{-1}g \in U$ for $g,h \in U'$. 
In view of Theorem~\ref{thm:kerana}, the analyticity of $\pi^v(g) = \pi(g)v$ 
on $xU'$, $x \in G$, is equivalent to the 
analyticity of the function 
$$ (g,h) \mapsto \la \pi(xg)v, \pi(xh)v \ra = \pi^{v,v}(h^{-1}g) $$ 
on $U' \times U'$, which follows from the analyticity of 
$\pi^{v,v}$ on $U$. 
\end{prf}

In view of the GNS construction, we obtain 
an analog of \cite[Cor.~7.4]{Ne10} for analytic 
functions: 
\begin{cor} If $\phi$ is a positive definite function on a Lie group 
$G$ which is analytic in a $\1$-neighborhood, then $\phi$ is analytic. 
\end{cor}

\begin{prob} Suppose that $(\pi, V)$ is a smooth representation 
of the Banach--Lie group $G$ on the Banach space $V$. 
Is the space $V^\omega$ of analytic vectors dense? 
This is true if $\dim G < \infty$ for any continuous 
Banach representation (\cite[Thm.~4]{Nel59}, \cite{Ga60}), 
but the proof is based on 
the integrated representation to $L^1(G)$ and mollifying with the 
heat kernel. All these techniques are not available for general infinite 
dimensional groups. However, some Lie groups carry so-called heat 
kernel measures which behave in many respects like Haar measure 
with a heat kernel density (cf. \cite{DG08}). 
\end{prob}

\section{Positive analytic functionals} \label{sec:6}

In this section we turn to the $*$-aspects of analytic functionals 
on $U(\g)$. Throughout $U_\C(\g) := U(\g_\C) \cong 
U(\g)_\C$ denotes the 
enveloping algebra of the complexified Lie algebra $\g_\C$. 
It coincides with the complexification of the real algebra $U(\g)$ 
and we have 
$\Hom_\R(U(\g),E) \cong \Hom_\C(U_\C(\g),E)$ for every complex 
vector space~$E$. 

\begin{defn}
We write $*$ for the unique antilinear antiautomorphism of 
$U_\C(\g)$ with $x^* = - x$ for $x \in \g$. Then 
$(U_\C(\g),*)$ is an involutive algebra and a linear functional
$\lambda \in U_\C(\g)^*$ is called {\it positive} if 
$\lambda(D^*D) \geq 0$ for $D \in U_\C(\g)$. 
\end{defn} 

Typical examples of positive functional arise from 
$*$-representations $(\rho, \cD)$ on pre-Hilbert spaces 
(see Definition~\ref{def:5.1}(b)(ii)) 
as $\lambda(D) = \la \rho(D)v,v\ra$, and a suitable variant 
of the GNS construction implies that every positive 
functional on $U_\C(\g)$ is of this from 
(\cite[Thm.~6.3]{Po71}, \cite[Thm.~8.6.2]{Sch90}). 

\begin{defn} 
A positive linear functional 
$\lambda \in U_\C(\g)^*$ is said to be {\it integrable} if 
there exists a smooth positive definite function 
$\phi \: G \to \C$ with 
$$ \lambda(D) = (D\phi)(\1) \quad \mbox{ for } \quad D \in U_\C(\g). $$

If $(\pi, \cH)$ is a continuous unitary representation 
and $v \in \cH$ with $\phi(g) = \la \pi(g)v,v\ra$ 
(the existence follows from the GNS construction), then 
\cite[Thm.~7.2]{Ne10} implies that $v$ is a smooth vector. 
Therefore the integrability of $\lambda$ is equivalent to the 
existence of a unitary representation $(\pi,\cH)$ and a 
smooth vector $v \in \cH^\infty$ with 
$$ \lambda(D) = \la \dd\pi(D)v,v \ra \quad \mbox{ for } \quad 
D \in U_\C(\g). $$
\end{defn}

The following proposition builds a bridge between 
analytic vectors and analytic functionals. 

\begin{prop}
  \label{prop:2.1a} 
Let $(\rho, \cD)$ be a strongly continuous $*$-representation of $\g$
on the pre-Hilbert space $\cD$ and  
$\rho \: U_\C(\g) \to \End(\cD)$ its canonical extension. 
For $v \in \cD$, the positive functional $\beta(D) := \la \rho(D)v,v\ra$ 
is analytic if and only if $v$ is an analytic vector. 
This implies that $v$ is an analytic vector for 
every operator $\rho(x)$, $x \in \g$, and that 
$\sum_{n = 0}^\infty \frac{\|\rho(x)^n v\|}{n!}< \infty$ 
for $p(x) < R_{\beta,p}/2$. 
\end{prop}

\begin{prf} First we observe that, for $x \in \g$, we have  
  \begin{equation}
    \label{eq:ident}
\|\rho(x)^nv\|^2 = \la \rho(x)^n v, \rho(x)^n v\ra = (-1)^n \beta(x^{2n}).   
\end{equation}
If $\beta$ is analytic, then Proposition~\ref{prop:anachar} 
provides an estimate of the form 
$\|\beta_n^s\|_p\leq C n! t^{-n}$ for every $t < R_{\beta,p}$. 
This leads to 
$$ \frac{\|\rho(x)^nv\|}{n!}
\leq \sqrt{C} \frac{\sqrt{(2n)!}}{n!} p(x)^n t^{-n}. $$
Therefore $\frac{\sqrt{(2n+2)(2n+1)}}{n+1} \to 2$
implies that 
$\sum_{n = 0}^\infty \frac{\|\rho(x)^nv\|}{n!}$
converges for $p(x) < t/2$, hence for $p(x) < R_{\beta,p}/2$. 

Suppose, conversely, that $\sum_{n = 0}^\infty \frac{\|\rho(x)^n v\|}{n!}$ 
converges for every $x$ in a $0$-neighborhood $U$ of $\g$. Then the 
series 
$$ \Big\la \sum_{n = 0}^\infty \frac{\rho(x)^n v}{n!} , v\Big\ra 
= \sum_{n = 0}^\infty \frac{\la \rho(x)^n v, v \ra}{n!}
= \sum_{n = 0}^\infty \frac{\beta(x^n)}{n!} $$
converges for every $x \in U$. 

To see that $v$ is an analytic vector for $\rho(x)$, pick 
$t > 0$ such that $tx \in U$ and observe that the series 
$$\sum_{n = 0}^\infty \frac{\|\rho(x)^nv\|t^n}{n!} = 
\sum_{n = 0}^\infty \frac{\|\rho(tx)^nv\|}{n!}$$
converges. 
\end{prf}

\begin{defn} (cf.\ \cite[Thm.~III.1.3]{Ne00}, \cite[Thm.~6.3]{Po71}, 
\cite[Thm.~8.6.2]{Sch90}) 
Let $\lambda \in U_\C(\g)^*$ be a positive analytic functional. 
On the dual space $U_\C(\g)^*$ we consider the {\it right regular representation} of $U_\C(\g)$ by $D\beta := \beta \circ \rho_D$. 
Then there exists a unique Hilbert space 
$\cH_\lambda \subeq U_\C(\g)^*$ containing $\lambda$ 
for which the subspace 
$$ \cD_\lambda := U_\C(\g) \lambda $$
is dense and on which the scalar product satisfies 
\begin{equation}
  \label{eq:evald}
\beta(D) = \la \beta, D^*\lambda \ra \quad \mbox{ for } \quad 
\beta \in \cH_\lambda, D \in U_\C(\g). 
\end{equation}
In particular, we have 
$$ \la D_1\lambda, D_2\lambda \ra = \lambda(D_2^*D_1). $$

Now $\rho_\lambda(D)\beta := D\beta := \beta \circ\rho_D$ 
defines a $*$-representation 
of $U_\C(\g)$ on $\cD_\lambda$ by unbounded operators, and for 
each $x \in \g$, the operator $\rho_\lambda(x)$ on 
$\cD_\lambda$ is skew-hermitian. 
\end{defn}

\begin{prop}
  \label{prop:2.1} 
If $\lambda$ is a positive analytic functional on 
$U_\C(\g)$, then the following assertions hold: 
\begin{description}
\item[\rm(a)] $\lambda \in \cD_\lambda$ 
is an analytic vector for every $\rho_\lambda(x)$, $x \in \g$. 
\item[\rm(b)] If $\lambda$ is $p$-analytic for a submultiplicative 
continuous seminorm $p$ on $\g$, then 
$\cD_\lambda = \rho_\lambda(U_\C(\g))\lambda$ is equianalytic. 
More precisely, the series 
$\sum_{n = 0}^\infty \frac{\rho_\lambda(x)^nv}{n!}$ converges 
for $p(x)  < R_{\lambda,p}/2$ and $v \in \cD_\lambda$. 
\end{description}
\end{prop}

\begin{prf} (a) In view of Proposition~\ref{prop:2.1a}, this follows from 
$\lambda(D) = \la \rho_\lambda(D)\lambda, \lambda\ra$ for 
the representation $(\rho_\lambda, \cD_\lambda)$. That this 
representation is strongly continuous follows from the continuity of 
the quadratic maps 
$x\mapsto \la \rho_\lambda(x)^2 \rho_\lambda(D)\lambda, \rho_\lambda(D)\lambda \ra 
= \lambda(D^* x^2 D)$
(Lemma~\ref{lem:4.5}). 

(b) Let $v := \rho_\lambda(D)\lambda$ for some 
$D \in U_\C(\g)$ and define 
$$\beta(D') := \la \rho_\lambda(D')v,v\ra 
= \lambda(D^*D'D) \quad \mbox{ for } \quad D' \in U_\C(\g). $$ 
Then Theorem~\ref{thm:1.5} implies that 
$\beta$ is an analytic functional whose radius of convergence satisfies
$R_{\beta,p} \geq r := R_{\lambda,p}$. 
Therefore Proposition~\ref{prop:2.1a} 
implies that $v$ is an analytic vector for every operator
$\rho_\lambda(x)$, $x \in \g$, and the corresponding exponential 
series converges for $p(x) < R_{\lambda,p}/2 \leq R_{\beta,p}/2$. 
\end{prf}

\begin{defn} \label{def:6.5} Let $(\rho, \cD)$ be a $*$-representation 
of $U_\C(\g)$ on the pre-Hilbert space $\cD$. We call 
a subset $E \subeq \cD$ {\it equianalytic} if 
there exists a $0$-neighborhood 
$U \subeq \g$ such that 
$$ \sum_{n = 0}^\infty \frac{\|\rho(x)^n v\|}{n!} < \infty 
\quad \mbox{ for } \quad v \in E, x \in U. $$
This implies in particular that each $v \in E$ is an analytic 
vector for every $\rho(x)$, $x \in \g$.
\end{defn}

\begin{lem} \label{lem:2.3a} 
Let $(\rho, \cD)$ be a $*$-representation of 
the BCH--Lie algebra $\g$ on the dense equianalytic 
subspace $\cD$ of $\cH$. 
Then the closures $\oline{\rho(x)}$, $x \in \g$, generate 
unitary one-parameter groups and there exists a $0$-neighborhood 
$U \subeq \g$ on which the BCH product $*$ is defined, and 
the map 
$\tilde\rho \: \g \to \U(\cH), x \mapsto e^{\oline{\rho(x)}}$
satisfies 
\begin{equation}
  \label{eq:multip}
\tilde\rho(x *y) = \tilde\rho(x) \tilde\rho(y) \quad \mbox{ for }\quad x,y \in U.
\end{equation}
\end{lem}

\begin{prf} Nelson's Theorem~\ref{thm:nelson} implies that 
the operators $\rho(x)$ are 
essentially skew-adjoint, so that their closures generate unitary 
one-parameter groups and $\tilde\rho$ is defined. 

Let $U_1 \subeq \g$ be an open balanced $0$-neighborhood with the 
property that the 
function 
$$ U_1 \to \cH, \quad x \mapsto \tilde\rho(x) = \sum_{n = 0}^\infty \frac{\rho(x)^n}{n!} $$ 
is analytic for every $v \in \cD$ (Definition~\ref{def:6.5}). 
Let $U_2 \subeq U_1$ be a smaller convex $0$-neighborhood for which the 
BCH-multiplication defines an analytic map 
$$U_2 \times U_2 \to U_1, \quad 
(x,y) \mapsto x + y + \frac{1}{2}[x,y] + \cdots.$$
Then we have for $v \in \cD$ two analytic functions 
$$ F_{1/2} \: U_2 \times U_2 \to \C, \quad 
F_1(x,y) = \la \tilde\rho(y)v, \tilde\rho(-x)v\ra, \quad 
F_2(x,y) = \la \tilde\rho(x*y)v, v\ra.  $$
To show that they coincide, it suffices to verify that 
their Taylor polynomials in $(0,0)$ agree. 
Let $\beta(D) := \la \rho(D)v,v\ra$. 
Then, for $F_1$, the homogeneous term of degree $(m,n)$ is given by 
$$ T_{m,n}(F_1)(x,y) 
= \frac{1}{m!n!} 
\la \rho_\lambda(y)^nv, \rho_\lambda(-x)^m v\ra 
= \frac{1}{m!n!} \la \rho_\lambda(x^m y^n)v, v \ra 
= \beta\Big(\frac{x^m y^n}{m!n!}\Big). $$ 
For $F_2$ the homogeneous term of degree $(m,n)$ has the form 
$$ \sum_{k \leq n+m} \frac{1}{k!} T_{m,n}(\la \rho_\lambda(x*y)^k v,v\ra)
=  \beta\Big(\sum_{k \leq n+m} \frac{1}{k!} T_{m,n}((x*y)^k)\Big). $$ 
From the identity 
$e^x e^y = e^{x*y}$
in the completed free associative algebra $\cA$ with two generators 
$x$ and $y$ we derive that 
$$ \frac{x^my^n}{m!n!} = \sum_{k \leq n+m} \frac{1}{k!} T_{m,n}((x*y)^k)$$ 
holds in $\cA$, and this implies the corresponding relation 
in $U(\g)$ (\cite[Ch.~2, \S 6]{Bou89}). We conclude that 
$$ \la \tilde\rho(x)\tilde\rho(y)v, v \ra 
=  \la \tilde\rho(x * y)v, v\ra 
\quad \mbox{ for } \quad x,y \in U_2, v \in \cD.$$ 
Now polarization leads to 
$$ \la \tilde\rho(x)\tilde\rho(y)v, w \ra 
=  \la \tilde\rho(x * y)v, w\ra 
\quad \mbox{ for } \quad v, w \in \cD.$$ 
Since $\cD$ is dense in 
$\cH$, we obtain \eqref{eq:multip} for 
$U := U_2$. 
\end{prf}

The following theorem is our first main result on integrability 
of Lie algebra representations based on a sufficient supply of 
analytic vectors. 

\begin{thm} \label{thm:2.8} {\rm(Integrability Theorem)} 
Let $G$ be a $1$-connected BCH--Lie group with Lie algebra $\g$ and 
$(\rho, \cD)$ be a $*$-representation of $\g$ on the dense subspace 
$\cD$ of $\cH$. Assume that either 
\begin{description}
\item[\rm(a)] $G$ is Banach and $\cD$ consists of analytic vectors, or 
\item[\rm(b)] $\cD$ is equianalytic. 
\end{description}
Then there exists a 
unique unitary representation $(\pi, \cH)$ of $G$  
on the corresponding Hilbert 
space $\cH$ with $\dd\pi(x)\res_{\cD} = \rho(x)$ for 
$x \in \g$. 
\end{thm}

\begin{prf} With Nelson's Theorem~\ref{thm:nelson} we see that 
the operators $\rho(x)$, $x \in \g$, are essentially skew-adjoint, 
so that their closures generate unitary one-parameter groups. This leads 
to a map 
$\tilde\rho \: \g \to \U(\cH), x \mapsto e^{\oline{\rho(x)}}.$
We claim that 
\begin{equation}
  \label{eq:mult2}
\tilde\rho(x*y) = \tilde\rho(x)\tilde\rho(y) 
\end{equation}
holds for $x,y$ in some open $0$-neighborhood $U \subeq \g$. 
It suffices to show that 
\begin{equation}
  \label{eq:mult2b}
\tilde\rho(x*y)v = \tilde\rho(x)\tilde\rho(y)v \quad \mbox{ for every} 
\quad v \in \cD. 
\end{equation}

(a) Suppose that $G$ is Banach. 
Let $\cH_v^0 := \rho(U_\C(\g))v$ and $\cH_v$ be the closure of this subspace 
in $\cH$. 
Then $\lambda(D) := \la \rho(D)v,v\ra$ is a positive analytic functional 
on $U_\C(\g)$ and the map 
$$ \Gamma \: \cH_v^0 \to U_\C(\g)^*, \quad 
\Gamma(w)(D) := \la \rho(D)w, v\ra $$
is equivariant with respect to the right regular representation on 
$U_\C(\g)^*$ 
and satisfies $\Gamma(\cH_v^0) = \rho_\lambda(U_\C(\g))\lambda = \cD_\lambda$. 
From $\la D_1 v, D_2 v \ra = \lambda(D_2^* D_1)$ 
we further derive that $\Gamma$ is isometric. 
Therefore Proposition~\ref{prop:2.1}(b) implies that 
$\cH_v^0$ is equianalytic. It therefore suffices to consider 
case (b) to prove \eqref{eq:mult2b}. 

(b) If $\cD$ is equianalytic, then Lemma~\ref{lem:2.3a} applies 
and provides an open $0$-neigh\-bor\-hood 
$U \subeq \g$ with \eqref{eq:mult2}. 
Now we derive from \cite[Ch.~3, \S 6, Lemma 1.1]{Bou89} 
(for Banach--Lie groups) and from \cite{GN10} 
(for BCH--Lie groups) that there exists a unique 
homomorphism $\pi \: G \to \U(\cH)$ such that 
$\pi(\exp_G x) = \tilde\rho(x)$
holds for all elements $x$ in some $0$-neighborhood of $\g$.
This relation implies in particular that, for each $x \in \g$, 
the unitary one-parameter group 
$\pi(\exp tx)$ is generated by the closure of the 
operator $\rho(x)$, so that $\pi$ is uniquely 
determined by the representation $\rho$ of $\g$ on $\cD$. 
\end{prf}

\begin{lem}
  \label{lem:2.3} 
If $\lambda$ is a positive analytic functional and 
$U\subeq \g$ a $0$-neighborhood, then 
$\tilde\rho_\lambda(U)\lambda$ is total in $\cH_\lambda$. 
\end{lem}

\begin{prf} Let $\cK$ denote the closure of the subspace generated 
by $\tilde\rho(U)\lambda$. From 
$$ \tilde\rho_\lambda(x)\lambda = \sum_{n = 0}^\infty \frac{1}{n!}\rho_\lambda(x)^n \lambda $$
for $x$ sufficiently close to $0$ we derive that 
$\rho_\lambda(x)^n \lambda
= \frac{d}{dt^n}\Big|_{t = 0} \tilde\rho_\lambda(tx)\lambda
 \in \cK.$
Now the assertion follows from the fact that 
$U_\C(\g)$ is spanned by the elements of the form 
$x^n$, $x \in \g$, $n \in \N_0$. 
\end{prf}

If $(\pi, \cH)$ is a unitary representation of the analytic Lie group 
$G$ and $v \in \cH^\omega$ an analytic vector, then the corresponding 
functional $\pi^{v,v}(D) := \la \dd\pi(D)v,v\ra$ on $U_\C(\g)$ 
is positive and analytic. The following theorem provides a converse 
which is known for finite dimensional groups 
(cf.\ \cite[Prop.~2.3]{Go69}). 
Its strength lies in the fact that the information about 
the whole representation $(\pi_\lambda, \cH_\lambda)$ is 
encoded completely on the infinitesimal level in the analytic 
functional $\lambda$. 

\begin{thm} \label{thm:2.7} {\rm(Integrability Theorem for 
analytic functionals)} 
Let $G$ be a $1$-connected BCH--Lie group with Lie algebra $\g$ and 
$\lambda \in U_\C(\g)^*$ be a positive analytic functional. 
Suppose that $U(\g)\lambda \subeq \Hom(U(\g),\C)_{\omega,p,r}$ 
for some $p \in \cP(\g)$ and $r> 0$, which is automatically 
satisfied if $\lambda \in \Hom(U(\g),\C)_{\omega,p}$ 
for some submultiplicative $p$. 
Then there exists a unique unitary representation 
$(\pi_\lambda, \cH_\lambda)$ of $G$ for which 
the dense subspace $\cD_\lambda$ consists of analytic vectors and 
$$ \dd\pi_\lambda(x)\res_{\cD_\lambda} = \rho_\lambda(x)
\quad \mbox{ for } \quad x \in \g. $$
\end{thm}

\begin{prf} If  $\lambda \in \Hom(U(\g),\C)_{\omega,p}$ for a 
submultiplicative seminorm $p$, then Theorem~\ref{thm:1.5} 
implies that $\cD_\lambda = U_\C(\g)\lambda 
\subeq \Hom(U(\g),\C)_{\omega,p,r}$ 
for some $r > 0$. We assume that this condition is satisfied. 
Then 
$$U_1 := \{ x \in \g \: p(x) < r/2\}$$ 
is an open convex $0$-neighborhood on which the  
series $\sum_{n = 0}^\infty \frac{1}{n!} \rho_\lambda(x)^nv$ 
converges (Proposition~\ref{prop:2.1a}). 
In particular, $\cD_\lambda$ is equianalytic 
and Theorem~\ref{thm:2.8}(b) applies. 
\end{prf} 

Specializing to Banach--Lie groups, the preceding theorem simplifies to: 
\begin{cor} \label{cor:2.7b} Let $G$ be a $1$-connected Banach--Lie 
group with Lie algebra $\g$ and 
$\lambda \in U_\C(\g)^*$ be a positive analytic functional. 
Then there exists a unique unitary representation 
$(\pi_\lambda, \cH_\lambda)$ of $G$ for which 
the dense subspace $\cD_\lambda$ consists of analytic vectors and 
$$ \dd\pi_\lambda(x)\res_{\cD_\lambda} = \rho_\lambda(x)
\quad \mbox{ for } \quad x \in \g. $$
\end{cor}

\subsection*{Integrability without a given group} 

Since not every Banach--Lie algebra is the Lie algebra 
of a Lie group (\cite{EK64}), it is natural to ask to which extent 
the assumption of the existence of $G$ in Theorem~\ref{thm:2.8}   
is needed to integrate the Lie algebra representation. 
The following theorem provides a corresponding variant, where no 
a priori given group is needed (cf.\ \cite{AN08} for a 
non-linear variant of this theorem). 

\begin{thm} \label{thm:2.8b} {\rm(Second Integrability Theorem)} 
Let $\g$ be a BCH--Lie algebra and 
$(\rho, \cD)$ be a faithful strongly continuous 
$*$-representation of $\g$ on the pre-Hilbert 
space $\cD$ by skew-adjoint operators. Assume that 
$\cD$ consists of analytic vectors, so that the closure of 
each operator  $\rho(x)$ generates a unitary one-parameter group 
on the completion $\cH$ of $\cD$. Assume that the subgroup 
$$ \Gamma := \{ z \in \z(\g) \: e^{\oline{\rho(z)}} = \1\} $$
of the center $\z(\g)$ of $\g$ 
is discrete. Then there exists a connected Lie group 
$G$ and a continuous unitary representation $(\pi, \cH)$ 
of $G$ on the corresponding Hilbert 
space $\cH$ with $\dd\pi(x)\res_{\cD} = \rho(x)$ for 
$x \in \g$. 
\end{thm}

\begin{prf} With the same argument as in the proof of 
Theorem~\ref{thm:2.8}, we obtain a map 
$\tilde\rho \: \g \to \U(\cH)$ satisfying 
$\tilde\rho(x*y) = \tilde\rho(x)\tilde\rho(y)$
for $x,y$ in some open balanced $0$-neighborhood $U \subeq \g$. 
Shrinking $U$, we may even assume that 
$$ \tilde\rho(x*y*z) = \tilde\rho(x)\tilde\rho(y)\tilde\rho(z) \quad \mbox{ for } \quad 
x,y,z \in U. $$
For $z = -x$ this leads to 
$$ \tilde\rho(e^{\ad x}y) = \tilde\rho(x*y*(-x)) 
= \tilde\rho(x)\tilde\rho(y)\tilde\rho(-x)= \tilde\rho(x)\tilde\rho(y)\tilde\rho(x)^{-1}, $$ 
and by passing to the infinitesimal generators, we obtain 
$$ \rho(e^{\ad x}y) = \tilde\rho(x)\rho(y)\tilde\rho(x)^{-1}. $$ 

We consider the subgroup $G := \la \tilde\rho(\g) \ra \subeq \U(\cH)$ 
generated by $\tilde\rho(\g)$. 
If $x \in U$ satisfies $\tilde\rho(x) = \1$, then the preceding 
argument, together with the injectivity of $\rho$ implies 
that $e^{\ad x} =\id_\g$. 
Since $\z(\g)$ is open in $\{ x \in \g \: e^{\ad x} = \id_\g\}$ 
(cf. the Adjoint Integrability Theorem in \cite{GN10}), 
we may shrink $U$ in such a way that, for 
$x \in U$, the relation $e^{\ad x} = \id_\g$ implies $x \in \z(\g)$. 
In view of the discreteness of $\Gamma$, we may further shrink 
$U$ so that 
$$ \tilde\rho^{-1}(\1) \cap U =  \{0\}. $$

If $V \subeq U$ is an open balanced $0$-neighborhood with 
$V * V \subeq U$, then $\tilde\rho\res_V$ is injective. 
In fact, if $\tilde\rho(x) = \tilde\rho(y)$ for $x,y \in V$, 
then 
$$ \tilde\rho((-x)*y) = \tilde\rho(-x) \tilde\rho(y) = \tilde\rho(x)^{-1} \tilde\rho(y) = \1 $$
leads to 
$(-x)*y \in \tilde\rho^{-1}(\1) \cap U = \{0\}, $
We conclude that 
$$x = x * 0 = x * ((-x) * y) = (x * (-x)) * y = 0 * y = y $$
and hence that $\tilde\rho\res_V$ is injective. 

This means that $\tilde\rho(V) \subeq G$ carries a manifold structure 
for which $\tilde\rho\res_V$ is a diffeomorphism, 
and that multiplication and inversion are smooth in an identity neighborhood in $V$. Since $G$ is generated by the image of every $0$-neighborhood in $V$, 
it carries a unique Lie group 
structure  for which $\tilde\rho$ defines a chart in a neighborhood of $0$ 
(\cite[Thm.~3.3.2, Rem.~3.3.3]{GN10}). 

The curves $\gamma_x \: \R \to G$ defined by 
$\gamma_x(t) := \tilde\rho(tx)$ define smooth $1$-parameter groups 
with $\gamma_x'(0) = x$, where we identify $\g$ via $T_0(\tilde\rho)$ with $T_\1(G)$. 
Therefore $\tilde\rho$ is an exponential function of $G$. 
In particular,  $G$ is locally exponential. 

For every vector $v \in \cD$, the map $x \mapsto \tilde\rho(x)v$ 
is analytic in a $0$-neighborhood, and since $\tilde\rho$ is a local 
diffeomorphism, the vector $v$ has a smooth orbit map. 
For the  unitary representation 
$(\pi,\cH)$ of $G$ given by the inclusion map 
$\pi \: G \to \U(\cH)$, this means that it is smooth, 
hence in particular continuous. 
From the construction it follows that, for every $x \in \g$, 
the generator $\oline{\dd\pi}(x)$ of the unitary one-parameter group  
$t \mapsto \pi(\exp_G(tx))$ restricts to $\rho(x)$ on~$\cD$. 
\end{prf}

\begin{rem} For a strongly continuous representation 
$(\rho,V)$ of a Banach--Lie algebra $\g$, the kernel is 
a closed ideal, so that $\g/\ker \rho$ also is a Banach--Lie algebra. 
In the situation of Theorem~\ref{thm:2.8b}, we may thus pass 
to an injective Lie algebra representation for which the 
theorem can be applied. 
\end{rem}

Here is an instructive example, showing that the group 
$\Gamma$ in Theorem~\ref{thm:2.8b} is not always discrete. 

\begin{ex} Let $(X,\fS,\mu)$ be a probability space. For a measurable 
function $f \: X \to \C$ we define 
$$ \|f\| := \sup\Big\{ \frac{\|f\|_n}{\root n \of {n!}} \: n \in \N\Big\} 
\quad \mbox{ for } \quad \|f\|_n := \Big(\int_X |f|^n\, d\mu\Big)^{1/n}. $$
Writing $\cM(X,\mu)$ for the space of measurable functions $f \:X \to \C$ 
modulo those vanishing on the complement of a zero set, 
we obtain a Banach space 
$$ V := \{ f \in \cM(X,\mu) \: \|f\| < \infty\}, $$
which can be viewed as a closed subspace of a weighted 
$\ell^\infty$-product of the spaces $L^n(X,\mu)$, $n \in \N$. 

For $\|f\| < \frac{1}{2}$ we have 
$$ \sum_n \frac{\|f^n\|_2}{n!} 
=  \sum_n \frac{\|f\|_{2n}^{n}}{n!} 
\leq  \sum_n \frac{\sqrt{(2n)!}}{n!} \|f\|^{n}
< \infty $$
because 
$\frac{\sqrt{(2n+2)(2n+1)}}{n+1} \to 2.$
Next we observe that for every measurable subset 
$E \subeq X$ we have $\chi_E \in V$ because 
$$ \|\chi_E\| 
= \sup_n \frac{\|\chi_E\|_n}{{\root n \of{n!}}} 
= \sup_n {\root n \of {\frac{\mu(E)}{n!}}} < \infty $$
follows from $\mu(E)/n! \to 0$. 
From 
$$ (n!)^{n+1} = (n!)^n n! \leq (n!)^n (n+1)^n 
= ((n+1)!)^n $$ 
we derive that ${\root n \of {n!}} \leq {\root {n+1} \of {(n+1)!}}$, 
i.e., the sequence ${\root n \of {n!}}$ is increasing. 
The formula for the radius of convergence of 
the exponential series now yields 
${\root n \of {n!}} \to \infty$.

We claim that $\chi_{E_k} \to 0$ in $V$ if 
$\mu(E_k) \to 0$. Let $\eps > 0$ and $N \in \N$ with 
${\root N \of{N!}} > \eps^{-1}$. Suppose that 
$\mu(E_k)^{1/n} \leq \eps$ for $k \geq N_0$ and every $n \leq N$. 
The relation 
$${\root n \of {\frac{\mu(E_k)}{n!}}} 
\leq \frac{1}{\root n \of {n!}} 
\quad \mbox{ for } \quad 
k \geq N_0, n > N$$
and 
$${\root n \of {\frac{\mu(E_k)}{n!}}} 
\leq \frac{\eps}{\root n \of {n!}} < \eps 
\quad \mbox{ for } \quad 
k \geq N_0, n \leq  N$$
then leads to 
$$ \|\chi_{E_k}\| \leq 
\sup_n {\root n \of {\frac{\mu(E_k)}{n!}}} \leq \eps 
\quad \mbox{ for } \quad k \geq N_0. $$
We conclude that $\chi(E_k) \to 0$ in $V$. 
Therefore the closed subgroup 
$$ \Gamma := \{ f \in V \: f(X) \subeq 2\pi \Z\} $$
is not discrete. 

The Banach space $V$ acts by unbounded multiplication operators 
on $L^2(X,\mu)$: 
$$ \rho(f)\xi := i f \xi. $$
The definition of the norm in $V$ implies that each 
bounded function $\xi \in L^2(X,\mu)$  is an analytic vector 
for $V$: 
$$ \sum_n \frac{\|(if)^n \xi\|_2}{n!} 
\leq \|\xi\|_\infty \sum_n \frac{\|f^n \|_2}{n!} < \infty 
\quad \mbox{ for } \quad \|f\| < \shalf. $$
Therefore we have a $*$-representation 
$(L^\infty(X,\mu), \rho)$ of $V$ by unbounded operators on the 
dense subspace $L^\infty(X,\mu)$ of the Hilbert space 
$L^2(X,\mu)$, for which the group 
$$ \{ f \in \g \: e^{\rho(f)} = \id\} = \Gamma $$ 
is not discrete. 
\end{ex}

\begin{rem}
The preceding example shows in particular that there exists an 
abelian Banach--Lie group $G$ and a continuous unitary representation 
$(\pi, \cH)$ for which $\cH^\omega$ is dense but 
the closed subgroup $\ker \pi$ is not discrete and satisfies 
$\L(\ker \pi) = \{0\}$. 
\end{rem}

\section{Extension of local positive definite functions}
\label{sec:7}

\begin{defn} Let $G$ be a group and 
$U \subeq G$ be a subset. A function $\phi \: UU^{-1} \to \C$ 
is said to be {\it positive definite} if the kernel 
$K \: U \times U \to \C, (x,y) \mapsto \phi(xy^{-1})$
is positive definite. 
\end{defn}

The main result of this section is Theorem~\ref{thm:extension} 
which asserts that any positive definite analytic function 
defined in a $\1$-neighborhood of a simply connected BCH--Lie group 
$G$ extends to an analytic function on $G$ which is positive 
definite. This is a quite remarkable 
result which is far from being true in the $C^\infty$-context. 
For a discussion of related extension problems for smooth 
positive definite functions on finite dimensional Lie groups 
we refer to \cite{Jo91} and the references given there. 
As we have seen in the introduction, for 
$G = \R^n$, such extension problems are part of the classical 
theory of moments. The corresponding 
version of Theorem~\ref{thm:extension} can be found in 
\cite[Prop.~3.10]{Fa08}. 

Before we turn to the extension theorem, we make the 
following observation on analytic functions on product 
spaces. 

\begin{lem} \label{lem:6.8} 
Let $X, Y$ be locally convex spaces and $Z$ a Banach space. 
Let $U \subeq X \times Y$ be an open neighborhood of $(x_0, y_0)$ 
and $f \: U \to Z$ be an analytic map. Then there exists an 
open neighborhood $U_{x_0}$ of $x_0$ and an open neighborhood 
$U_{y_0}$ of $y_0$ such that for every $x \in U_{x_0}$ 
and $h \in U_{y_0} - y_0$ we have 
$$ f(x, y_0 + h) = \sum_{n =0}^\infty \frac{1}{n!} 
(\partial_{(0,h)}^n f)(x,y_0). $$
\end{lem}

\begin{prf} In view of \cite[Thm.~5.1,Prop.~5.4]{BS71b}, we may 
assume that 
$f$ has a complex analytic extension 
$f \: U_\C \to Z_\C$, where 
$U_\C$ is an open neighborhood of $(x_0, y_0)$ in the complex locally 
convex space 
$X_\C \times Y_\C$. Let $V_{x_0}$, resp., $V_{y_0}$ 
be an open balanced neighborhood of $x_0$ in $U_\C$, resp., 
$y_0$ in $Y_\C$ such that 
$V_{x_0} \times V_{y_0} \subeq U_\C$. Then 
\cite[Prop.~5.5]{BS71b} implies that 
$$ f(x, y_0 + h) = \sum_{n =0}^\infty \frac{1}{n!} 
(\partial_{(0,h)}^n f)(x,y_0) 
\quad \mbox{ for } \quad x \in V_{x_0}, h \in V_{y_0} - y_0.$$ 
Hence the assertion follows with 
$U_{x_0} := V_{x_0} \cap X$ and 
$U_{y_0} := V_{y_0} \cap Y$.
\end{prf}

\begin{thm} {\rm(Extension of local positive definite analytic functions)} 
\label{thm:extension} 
Let $G$ be a simply connected Fr\'echet--BCH--Lie group, 
$V \subeq G$ an open connected $\1$-neighborhood and 
$\phi \: VV^{-1} \to \C$ an analytic positive definite function.  
Then there exists a unique positive definite 
analytic function $\tilde\phi \: G \to \C$ extending $\phi$. 
\end{thm}

\begin{prf} The uniqueness of $\tilde\phi$ follows from the connectedness 
of $G$ and the uniqueness of analytic continuation. 

{\bf Step 1:} To obtain its existence, we consider the reproducing kernel 
Hilbert space $\cH_K \subeq \C^V$ defined by the kernel $K$ via 
$f(g) = \la f, K_g \ra$ for $g \in V$ and $K_g(h) = K(h,g)$. 
Then the analyticity of the function 
$$ V \times V \to \C, \quad (g,h) \mapsto 
\la K_h, K_g \ra = K(g,h) = \phi(gh^{-1}) $$
implies that the map $\eta \: V \to \cH_K, g \mapsto K_g$ 
is analytic (Theorem~\ref{thm:kerana}). Here we use that 
$G$ is assumed to be Fr\'echet. Hence all functions in
 $\cH_K$ are analytic, and that we obtain for each 
$x \in \g$ an operator on $\cH_K^0 := \Spann \eta(V)$ by 
$$ (\rho(x)f)(g) := \derat0 f(g\exp_G(tx)) $$
because $\frac{d}{dt}|_{t=0} K_y(g\exp_G(tx))$ exists for every $g \in V$ and $x \in \g$. 

{\bf Step 2:} That we thus obtain a representation 
$\rho \: \g \to \End(\cH_K^0)$ follows from the embedding 
$\cH_K^0 \into C^\omega(V,\C)$ and the fact that 
$\g$ acts by a Lie algebra of (left invariant) 
vector fields on this space. Next 
we observe that 
\begin{align*}
& \la \rho(x)K_g, K_h \ra 
= \derat0 K_g(h\exp_G(tx)) 
= \derat0 K(h\exp_G(tx), g) \\
&= \derat0 \phi(h \exp_G(tx) g^{-1})
= \derat0 K(h, g\exp_G(-tx)) 
= \derat0 \oline{K_h(g\exp_G(-tx))}\\
&= -  \oline{(\rho(x)K_h)(g)} 
= - \oline{\la \rho(x)K_h, K_g \ra} 
= - \la K_g, \rho(x)K_h\ra. 
\end{align*}
Therefore $\rho$ extends to a $*$-representation of $U_\C(\g)$ on 
$\cH_K^0$. 

{\bf Step 3:} We claim that each $\eta(g) = K_g$, $g \in V$, 
is an analytic vector with 
\begin{equation}
  \label{eq:step3}
\eta(g \exp_G(-x)) = \sum_{n = 0}^\infty \frac{1}{n!} 
\rho(x)^n \eta(g) 
\end{equation}
for $x$ sufficiently close to $0$. 
The function $x \mapsto \eta(g\exp_G x)$ is analytic on a $0$-neighborhood 
of $\g$ and 
\begin{align*}
\frac{d^n}{dt^n}\Big|_{t=0} \big(\eta(g\exp_G(tx))\big)(y) 
&=  \frac{d^n}{dt^n}\Big|_{t=0} K_{g\exp_G tx}(y) 
=  \frac{d^n}{dt^n}\Big|_{t=0} \phi(y\exp_G(-tx)g^{-1})\\
&=  \frac{d^n}{dt^n}\Big|_{t=0} K_{g}(y \exp_G(-tx)) 
=  (\rho(-x)^n \eta(g))(y). 
\end{align*}
This leads to 
$$ \frac{d^n}{dt^n}\Big|_{t=0} \eta(g\exp_G(-tx)) 
= \rho(x)^n \eta(g), $$
and hence that \eqref{eq:step3} holds on a $0$-neighborhood of~$\g$. 

{\bf Step 4:} Let $W_G \subeq V$ be an open $\1$-neighborhood 
and $W_\g \subeq \g$ an open balanced $0$-neighborhood 
with $W_G \exp_G(W_\g) \subeq V$. Next we show that 
$\eta(W_G)$ is equianalytic and spans a dense subspace of $\cH_K$. 

Since the map $W_G \times W_\g\to \cH_K, (g,x) \mapsto 
\eta(g \exp_G(x))$ is analytic, Lemma~\ref{lem:6.8} shows that, 
after shrinking $W_G$ and $W_\g$, we may assume that 
$$ \eta(g \exp_G (-x)) = \sum_{n = 0}^\infty 
\frac{1}{n!} \rho(x)^n \eta(g) 
\quad \mbox{ for } \quad g \in W_G, x \in W_\g. $$
This means that $\eta(W_G)$ is equianalytic. 
To see that $\cK := \oline{\Spann(\eta(W_G))}$ is dense, 
we use the analyticity of $\eta$ to see that 
$\eta(V) \subeq \cK$, which shows that $\cH_K = \cK$. 

Applying Theorem~\ref{thm:2.8} 
with $\cD = \Spann(\eta(W_G))$, we now obtain a 
continuous unitary representation 
$(\pi, \cH_K)$ with 
$\pi(\exp_G(x)) = e^{\oline{\rho(x)}}$ for every $x \in \g$. 
Then 
$$\tilde\phi(g) := \la \pi(g)K_\1, K_\1\ra
= \la \pi(g) \phi, \phi \ra $$ 
is an analytic positive definite function on $G$, 
and for $x \in W_\g$ we have 
$$ \tilde\phi(\exp_G x) 
= \sum_{n = 0}^\infty \frac{\la \rho(x)^n \phi, \phi\ra}{n!} 
= \la \eta(\exp_G(-x)), \phi \ra 
= K_{\exp_G(-x)}(\1) = \phi(\exp_G x). $$
Therefore 
$\tilde\phi$ coincides with $\phi$ in a $\1$-neighborhood, so that 
the analyticity of $\phi$ and the connectedness of $V$ 
lead to $\tilde\phi\res_V = \phi$.
\end{prf}

{\bf Acknowledgement:} We thank St\'ephane Merigon for a careful 
reading of this paper and for several discussions on its subject matter. 
We are also grateful to J.~Faraut for pointing out some crucial 
references.

\end{document}